\documentclass[draft,preprint,12pt,numbers,sort&compress]{elsarticle}
\usepackage{}
\usepackage{mathrsfs}
\usepackage{amssymb}
\usepackage{amsthm}
\usepackage{mathrsfs}
\usepackage[centertags]{amsmath}
\usepackage{amsfonts}

\usepackage{color}

\usepackage{ctex}
\usepackage{CJK}

\input amssym.def
\input amssym.tex

\date{}

\newtheorem{Theorem}{Theorem}[section]
\newtheorem{Definition}{Definition}[section]
\newtheorem{Proposition}{Proposition}[section]
\newtheorem{Corollary}{Corollary}[section]
\newtheorem{Lemma}{Lemma}[section]

\newcommand\R{\mbox{\bf R}}

\newcommand\SR{\mbox{\scriptsize\bf R}}

\newcommand{\definition}{{\lower .5ex
  \hbox{$\>\>\stackrel{\triangle}{=}\>\>$} }}

\newcommand {\dps} {\displaystyle}

\begin{document}

\baselineskip=22pt
\thispagestyle{empty}

\begin{center}
{\Large \bf Sharp well-posedness of the Cauchy problem for}\\[1ex]
{\Large \bf the rotation-modified Kadomtsev-Petviashvili}\\[1ex]
{\Large \bf equation in anisotropic Sobolev spaces}\\[3ex]

{Wei Yan\footnote{Emails:   011133@htu.edu.cn; zhangyimin@whut.edu.cn;
yshli@scut.edu.cn; duan@iit.edu.}$^{a}$, \quad
Yimin Zhang$^{b}$,\quad
Yongsheng Li$^{c}$,\quad
Jinqiao Duan$^{d*}$}\\[1ex]

{$^a$School of Mathematics and Information Science, Henan
Normal University,}\\
{Xinxiang, Henan 453007,   China}\\[1ex]

{$^b$Department of Mathematics, Wuhan University of  Technology, }\\
{Wuhan, Hubei 430070,   China}\\[1ex]

{$^c$School of Mathematics, South  China  University of  Technology,}\\
{Guangzhou, Guangdong 510640,    China}\\[1ex]

{$^d$Department of Applied Mathematics, Illinois Institute
 of Technology,}\\[1ex]
{Chicago, IL 60616, USA  }\\[1ex]

\end{center}
\noindent{\bf Abstract.} We consider the Cauchy problem for
 the rotation-modified
 Kadomtsev-Petviashvili (RMKP) equation
\begin{align*}
 \partial_{x}\left(u_{t}-\beta\partial_{x}^{3}u
+\partial_{x}(u^{2})\right)+\partial_{y}^{2}u-\gamma u=0
 \end{align*}
  in the anisotropic Sobolev spaces $H^{s_{1},\>s_{2}}(\R^{2})$.
   When $\beta <0$ and  $\gamma >0,$  we prove that
the Cauchy problem is locally well-posed in
$H^{s_{1},\>s_{2}}(\R^{2})$ with $s_{1}>-\frac{1}{2}$
and $s_{2}\geq 0$. Our result considerably improves the Theorem 1.4 of
 R.  M. Chen, Y. Liu, P. Z. Zhang(
 Transactions of  the American Mathematical Society, 364(2012), 3395--3425.).
The key idea is that we divide the frequency
space into regular region and singular region.
We further  prove that  the Cauchy problem for RMKP equation is
 ill-posed in $H^{s_{1},\>0}(\R^{2})$
  with $s_{1}<-\frac{1}{2}$ in the sense that the
 flow map associated to the rotation-modified Kadomtsev-Petviashvili
 is not $C^{3}$.
 When $\beta <0,\gamma >0,$ by using the $U^{p}$ and $V^{p}$ spaces,
   we prove that the Cauchy problem is locally well-posed in
$H^{-\frac{1}{2},\>0}(\R^{2})$.

 \bigskip


\leftskip 0 true cm \rightskip 0 true cm

\newpage

\baselineskip=20pt

\bigskip
\bigskip

{\large\bf 1. Introduction}
\bigskip

\setcounter{Theorem}{0} \setcounter{Lemma}{0}\setcounter{Definition}{0}\setcounter{Proposition}{1}

\setcounter{section}{1}

In this paper, we investigate the Cauchy problem for  the rotation-modified
   Kadomtsev-Petviashvili (RMKP)   equation
\begin{align}
& \partial_{x}\left(u_{t}-\beta\partial_{x}^{3}u+\partial_{x}(u^{2})\right)
+\partial_{y}^{2}u-\gamma u=0,\label{1.01}\\
&u(x,y,0)=u_{0}(x,y),   \label{1.02}
\end{align}
in an anisotropic Sobolev space $H^{s_{1},s_{2}}(\R^{2})$.
This  equation is a model
 describing small amplitude, long internal waves in a rotating
 fluid propagating in one dominant direction with slow transverse effects,
  where the effects of rotation balance with weakly nonlinear
 and dispersive effects \cite{G1985,G1998}. Here $u(x,y,t)$ can be
regarded as the free surface,
 $t\in \R^{+}$ is a time-like variable, $x\in \R$,
is a spatial variable in the dominant direction of wave propagation and $y\in \R$ is a spatial variable
 in a direction transverse to the $x$-direction. The coefficient $\beta$
 determines the type of dispersion; in the case $\beta<0$
(negative dispersion),
 the equation models gravity surface waves in a shallow water channel
and internal
  waves in the ocean, while in the case $\beta>0$ (positive dispersion)
it models capillary surface waves or oblique magnet-acoustic waves in plasma.
The parameter $\gamma>0$ measures the effects of rotation and is proportional
to the Coriolis force.

When $\gamma=0$, (\ref{1.01}) reduces to the  Kadomtsev-Petviashvili
 (KP)   equation
\begin{align}
   \partial_{x}\left(u_{t}
  -\beta\partial_{x}^{3}u
  +\partial_{x}(u^{2})\right)
  +\partial_{y}^{2}u=0.\label{1.03}
\end{align}
When $\beta<0$, (\ref{1.03}) is the KP-I equation, while when $\beta>0$,
it  is the KP-II equation.
The KP-I and KP-II equations arise in physical contexts as models for
 the propagation
of dispersive long waves with weak transverse
effects \cite{KP, AS,KB1991,KB-1991},    and they are two-dimensional
 extensions
 of the Korteweg-de-Vries
equation \cite{KP}.

  Many authors have investigated the Cauchy problem for KP equations,
  \cite{Bourgain-GAFA-KP,BS,CKS-GAFA,CIKS,GPW,Hadac2008,Hadac2009,
 HN,IMCPDE,IMEJDE,IKT,ILP,IMS,IMT,IM2011,Kenig,KL,MST-Duke,MST2002,
 MST2004,MST2007,MST2011,TDCDS,TADE,TT,TCPDE,TDIE,TIMRN,Z}.
  Bourgain \cite{Bourgain-GAFA-KP} established
   the global well-posedness of  the Cauchy problem for the
 KP-II equation  in $L^{2}(\R^{2})$ and $L^{2}(\mathbf{T}^{2}).$
 Takaoka and  Tzvetkov \cite{TT} and
 Isaza and  Mej\'{\i}a \cite{IMCPDE} established the local well-posedness
  of KP-II equation in $H^{s_{1},s_{2}}(\R^{2})$
 with $s_{1}>-\frac13$ and $s_{2}\geq0.$ Takaoka \cite{TDCDS}
established the local well-posedness
 of KP-II equation in $H^{s_{1},0}(\R^{2})$
 with $s_{1}>-\frac12$ under the assumption that
$D_{x}^{-\frac{1}{2}+\epsilon}u_{0}\in L^{2}$
 with a suitably chosen $\epsilon$, where $D_{x}^{-\frac{1}{2}+\epsilon}$
 is Fourier multiplier
 operator with multiplier $|\xi|^{-\frac{1}{2}+\epsilon}.$ By using the $U^{p}$ and
 $V^{p}$ spaces,
 Hadac  et  al. \cite{Hadac2009}  established the small data
 global well-posedness
 and scattering result of  KP-II equation in the homogeneous
 anisotropic Sobolev space
$\dot{H}^{-\frac{1}{2},\>0}(\R^{2})$,  and arbitrary large
 initial data local well-posedness
 in both homogeneous Sobolev space $\dot{H}^{-\frac{1}{2},\>0}(\R^{2})$ and
  inhomogeneous anisotropic  Sobolev space
$H^{-\frac{1}{2},\>0}(\R^{2})$.

The resonant function of KP-I equation does not possess the
 good property as that of KP-II equation.
For the KP-I equation, it is  known that the Cauchy problem
 is globally well-posed
in the  second energy spaces
\begin{equation*}
    E^{2}=\left\{u_{0}\in L^{2}(\R^{2}): \|u_{0}\|_{L^{2}}+\|\partial_{x}u_{0}\|_{L^{2}}+\|\partial_{x}^{2}u_{0}\|_{L^{2}}
    +\|\partial_{x}^{-1}\partial_{y}u_{0}\|_{L^{2}}
    +\left\|\partial_{x}^{-2}\partial_{y}^{2}u_{0}\right\|_{L^{2}}<\infty\right\}
\end{equation*}
on both $\R^{2}$ and $\mathbf{T^{2}}$
 \cite{Kenig,MST2002,MST2004}.
For KP-I equation, Molinet et  al. \cite{MST-Duke} proved that the Picard
  iteration methods fails in standard Sobolev space and in anisotropic
Sobolev space,  since the flow map fails to be real-analytic at the
 origin in these spaces.
By  introducing  some  resolution  spaces and bootstrap  inequality
as well as the energy estimates,
Ionescu et al. \cite{IKT} established the  global  well-posedness of
KP-I in the natural energy space
$$E^{1}=\left\{u_{0} \in L^{2}(\R^{2}),\quad\partial_{x}u_{0}
\in L^{2}(\R^{2}),
\quad\partial_{x}^{-1}\partial_{y}u_{0}\in L^{2}(\R^{2})\right\}.$$
Molinet et al. \cite{MST2007} proved that the Cauchy problem for
 the KP-I equation
 is locally well-posed
in $H^{s,\>0}(\R^{2})$ with $s>\frac{3}{2}.$
Guo et al. \cite{GPW} proved that it is also true for $s=1.$
  Zhang \cite{Z} considered that periodic  KP-I equation and established
 local well-posedness
in the Besov type space
$
B_{2,1}^{\frac{1}{2}}(\mathbf{T}^{2}).
$

The existence and stability of solitary wave solutions of
 (\ref{1.01}) with $\beta>0$ were considered by Chen et al.\cite{CHL}.
Esfahani and  Levandosky \cite{EL} studied the strong instability of
\begin{align*}
&& \partial_{x}\left(u_{t}-\beta\partial_{x}^{3}u+\partial_{x}(u^{p})\right)
+\partial_{y}^{2}u=0,\quad p\geq2.
\end{align*}
Chen et al. \cite{CLZ} studied the solitary waves of the RMKP
equation  and   proved that
the Cauchy problem for  the rotation-modified   Kadomtsev-Petviashvili
 equation is
locally well-posed in the anisotropic Sobolev space $H^{s_{1},s_{2}}(\R^{2})$
with $s_{1}>-\frac{3}{10},s_{2}\geq0$, $\beta<0,\gamma>0$ and
 globally well-posed in the space $L^{2}(\R^{2})$.

 In this paper, we are going to study (\ref{1.01}) with $\beta<0$ and
 $\gamma>0$.
 Applying the operator $\partial_{x}^{-1}$ to the both sides of (\ref{1.01})
  yields
\begin{align}
u_{t}-\beta\partial_{x}^{3}u+\partial_{x}(u^{2})+
\partial_{x}^{-1}\partial_{y}^{2}u
-\gamma \partial_{x}^{-1}u=0.\label{1.04}
\end{align}
Applying the scaling transformation $v(x,y,t)=\beta^{-1}
u(x,(-\beta)^{-\frac{1}{2}}y,-\beta^{-1}t)$,  without loss of generality,
we may assume
 that $\gamma=1=-\beta.$
Motivated by   \cite{Bourgain93, Hadac2008,KPV1993,KPV1996},
we are going to prove that the Cauchy problem for (\ref{1.04})
is locally well-posed  in the anisotropic Sobolev spaces
$H^{s_{1},\>s_{2}}(\R^{2})$
  with $s_{1}>-\frac{1}{2}$ and $s_{2}\geq 0$.
   Thus, our result considerably improves the Theorem 1.4 of \cite{CLZ}.
  Motivated by Theorem 4.2 of \cite{KZ}, we shall prove also that
 (\ref{1.04}) is ill-posed in the anisotropic Sobolev
 spaces $H^{s_{1},\>0}(\R^{2})$
  with $s_{1}<-\frac{1}{2}$  in the sense that
   the flow map
 associated to the RMKP is not $C^{3}$.
 Inspired by \cite{Hadac2009}, when $\beta<0,\gamma>0,$
 by using the $U^{p}$ and $V^{p}$ spaces,  we prove that
 the Cauchy problem is locally well-posed in
$H^{-\frac{1}{2},\>0}(\R^{2})$.

\bigskip

Now we briefly discuss  the difficulties  in showing that
 the local well-posedness of  (\ref{1.04})
  in the anisotropic Sobolev spaces $H^{s_{1},\>s_{2}}(\R^{2})$
  with $s_{1}>-\frac{1}{2}$ and $s_{2}\geq 0$.
We note that the structure of (\ref{1.04}) is more complicated than that
 of the KP-II equation due to the presence of the rotation term.
 More precisely, the phase function of  KP-II equation is
$\xi^{3}-\frac{\eta^{2}}{\xi}$, while the phase function of the
RMKP-II equation is $\xi^{3}-\frac{\eta^{2}+1}{\xi}$.
In \cite{Hadac2008}, Hadac  introduced some suitable  function spaces
and establishing  Theorem 3.3  which plays the key role
in showing the sharp local well-posedness of KP-II equation,
  The integral in Theorem 3.3 was estimated in the whole frequency space.
  Due to the presence of the rotation term, i.e. the presence
 of an additional term $\frac{1}{\xi}$ in the phase function,
 the estimates for  the integral in  Theorem 3.3 in \cite{Hadac2008}
\begin{align}
   \left|\int_{\SR^{6}}
   \frac{
     |\xi_{1}|^{-\frac{1}{2}}
     |\xi_{2}|^{\frac{1}{2}}
     }{
     \langle\lambda_{1}\rangle^{b} \langle\lambda_{2}\rangle^{b}
     }
     F(\zeta)
     F_{1}(\zeta_1)
     F_{2}(\zeta_2)
    \,d\zeta_{1}d\zeta\right|
\leq C\|F\|_{L_{\zeta}^{2}}\>
      \|F_{1}\|_{L_{\zeta}^{2}}\>
      \|F_{2}\|_{L_{\zeta}^{2}}
      \label{1.05}
\end{align}
are no longer valid in the whole frequency space for (\ref{1.04}),
and thus we can not completely apply the method to our case.

To overcome this difficulty we have to split the  frequency space as follows.
In this paper we always write $
\zeta_1=(\xi_1,\eta_1,\tau_1),\, \zeta_2=(\xi_2,\eta_2,\tau_2),\,
\zeta =(\xi,\eta,\tau)=\zeta_1+\zeta_2$.
Let
\begin{align}
\Omega_0& =\bigl\{(\zeta_1,\zeta)\in
     \R^{6}\,\Bigl |\,
   3|\xi_{1}|\leq |\xi_{2}|,
    |\xi_{2}|\geq 6^{7}\bigr\}   ,\label{1.06}\\
\Omega_1&=\left\{(\zeta_1,\zeta)\in \R^6
     \,\biggl |\>
     \Big|1-\frac{1}{3\xi_{1}^2\xi_{2}^{2}}\Big|
       \geq \frac{1}{4}\right\},\label{Omega1}\\
\Omega_2 &=\left\{(\zeta_1,\zeta)\in \R^6
     \,\biggl |\>
   \left|1-\frac{1}{3\xi_{1}^2\xi_{2}^{2}}
   \right|< \frac{1}{4}\right\},\label{Omega2}\\
A&= \bigg\{(\zeta_1,\zeta)\in \R^6
    \,\bigg|\, |3\xi\xi_{1}\xi_{2}|
    \geq \frac{\xi_{1}^{2}-\xi_{1}\xi_{2}+\xi_{2}^{2}}
              {|\xi\xi_{1}\xi_{2}|}\bigg\},\label{RegionA}\\
B&= \bigg\{(\zeta_1,\zeta)\in \R^6
    \,\bigg|\, |3\xi\xi_{1}\xi_{2}|
    \leq \frac{\xi_{1}^{2}-\xi_{1}\xi_{2}+\xi_{2}^{2}}
              {|\xi\xi_{1}\xi_{2}|}\bigg\}.
              \label{RegionB}
\end{align}
We refer $\Omega_r=\Omega\bigcap\Omega_1$
the  regular region,  and $\Omega_s=\Omega_0\bigcap\Omega_2$
the singular region.
For the regular region $\Omega_r$, we can adopt the idea of \cite{Hadac2008}
to obtain Lemmas 2.4--2.6 in the present paper, which, together with
 the algebraic identity
\begin{align}
&\left|\xi_{1}^{3}-\frac{\eta_{1}^{2}+1}{\xi_{1}}+\xi_{2}^{3}-
\frac{\eta_{2}^{2}+1}{\xi_{2}}-(\xi^{3}-\frac{\eta^{2}+1}{\xi})\right|
\nonumber\\
=&\left|3\xi\xi_{1}\xi_{2}+\frac{(\xi_{1}\eta-\xi\eta_{1})^{2}}
{\xi\xi_{1}\xi_{2}}+\frac{\xi^{2}-\xi\xi_{1}+\xi_{1}^{2}}
{\xi\xi_{1}\xi_{2}}\right|,\label{1.09}
\end{align}
implies the bilinear estimates for the regular region part.
For the singular region part $\Omega_s$, we have
\begin{align}
|\xi_{1}|\sim |\xi_{2}|^{-1}.\label{1.010}
\end{align}
(\ref{1.05})
does not hold in $\Omega_s$, i.e. in the region of
the interaction between the very low frequency and the very high
frequency.  To bypass this problem we  establish the estimate
in $\Omega_{s}$ according to the different intervals of symbol function
  $\lambda=\tau-\phi(\xi,\eta)$, see Lemmas  2.9, 2.11 and 2.13 for details.

We now introduce some notations. Throughout this paper, we assume that
$C$ is a generic positive constant and may vary from line to line.
The notation  $a\sim b$ means that there exists constants $C_{j}>0(j=1,2)$
such that $C_{1}|b|\leq |a|\leq C_{2}|b|$, and $a\gg b$ means that
there exists a positive constant $C^{\prime}$ such that $|a|> C^{\prime}|b|.$
$0<\epsilon\ll1$ means that $0<\epsilon<\frac{1}{100}$.
Denote the phase function $\phi(\xi,\eta)=\xi^{3}-\frac{\eta^{2}+1}{\xi},$
the symbol function $\lambda=\tau-\phi(\xi,\eta)$.
Denote $\mathscr{F}$ the Fourier transformation in
$(x,y,t)$ and $\mathscr{F}^{-1}$ its inverse.
Denote $\mathscr{F}_{xy}$  the Fourier transformation in
$(x,y)$ and $\mathscr{F}^{-1}_{xy}$ its inverse. In the whole paper we denote
  \begin{align*}
  &\langle\cdot\rangle:=1+|\cdot|,\quad
  \lambda_{j} =\tau_{j}-\phi(\xi_{j}, \eta_{j})\ \ (j=1,2),\\
  &\zeta=(\xi,\eta,\tau),\quad
   \zeta_j=(\xi_j,\eta_j,\tau_j) \ \ (j=1,2),\\
  &d\zeta=d\xi d\eta d\tau ,\quad
   d\zeta_j=d\xi_j d\eta_j d\tau_j\ \ (j=1,2),\\
  &D_{x}^{a}u(x,y,t):=\frac{1}{(2\pi)^{\frac{3}{2}}}\int_{\SR^{2}}|\xi|^{a}
  \mathscr{F}u(\zeta)e^{ix\xi+iy\eta+it\tau} d\zeta ,\\
 & W(t)f:=\frac{1}{2\pi}\int_{\SR^{2}}e^{ix\xi+iy\eta+it\phi(\xi,\eta)}
 \mathscr{F}_{xy}f(\xi,\eta)d\xi d\eta,\nonumber\\
  &P^{N}u(x,y,t)=\frac{1}{(2\pi)^{\frac{3}{2}}}\int_{|\xi|\geq N}
  \mathscr{F}u(\zeta)e^{ix\xi+iy\eta+it\tau} d\zeta ,
\nonumber\\
  &P_{N}u(x,y,t)=\frac{1}{(2\pi)^{\frac{3}{2}}}\int_{|\xi|\leq N}
  \mathscr{F}u(\zeta)e^{ix\xi+iy\eta+it\tau} d\zeta ,
\nonumber
  \end{align*}
$\psi(t)$ is a smooth function supported in $[-2,2]$  and $\psi(t)\equiv 1$
in $[-1,1]$.  For a set $E$ we denote by ${\rm mes}(E)$ its
 Lebesgue measure and
$\chi_E$ its characteristic function. We will use the convention that capital letters denote
dyadic numbers, e.g. $N=2^{n}$  for $n\in Z$ and for a dyadic summation we write
$\sum\limits_{N}a_{N}=\sum\limits_{n\in Z}a_{2^{n}}$  and $\sum\limits_{N\geq M}a_{N}:=\sum\limits_{n\in Z:2^{n}\geq M}a_{2^{n}}$ for brevity. Let $\chi \in C_{0}^{\infty}((-2,2))$ be an even, non-negative function such that
$\chi(t)=1$ for $|t|\leq1$. We define $\Psi(t)=\chi(t)-\chi(2t)$ and $\Psi_{N}:=\Psi(N^{-1}\cdot)$. Then,
$\sum\limits_{N}\Psi_{N}(t)=1$ for $t\neq0$. We define $\mathscr{F}Q_{N}u:=\Psi_{N}\mathscr{F}u$
and
$\mathscr{F}Q_{0}u=\chi(2\cdot)\mathscr{F}u,Q_{\geq M}=\sum\limits_{N\geq M}Q_{N}$ as well as $Q_{<M}=I-Q_{\geq M}.$

 For $\xi,\xi_1,\xi_2\in \R$, we write
$$
|\xi_{\rm min}|:={\rm min}\left\{|\xi|,|\xi_{1}|,|\xi_{2}|\right\}
\quad \mbox{\rm and} \quad
 |\xi_{\rm max}|:={\rm max}\left\{|\xi|,|\xi_{1}|,|\xi_{2}|\right\}.
$$
$L_{t}^{r}L_{xy}^{p}$ is the function space with norm
 \begin{align*}
 \|f\|_{L_{t}^{r}L_{xy}^{p}}:=\left(\int_{\SR}\left(\int_{\SR^{2}}
 |f|^{p}dxdy\right)^{\frac{r}{p}}dt\right)^{\frac{1}{r}}.
 \end{align*}
 The anisotropic Sobolev space $H^{s_{1},s_{2}}(\R^{2})$ is
defined via its norm
 \begin{align*}
 \|u_{0}\|_{H^{s_{1},s_{2}}(\SR^{2})}=\left\|\langle\xi\rangle^{s_{1}}
 \langle\eta\rangle^{s_{2}}
 \mathscr{F}_{xy}u_{0}(\xi,\eta)\right\|_{L_{\xi\eta}^{2}}.
 \end{align*}
The space
$
  X_{\sigma}^{s_{1},\>s_{2},\>b}
$ is defined via its norm
$$
 \|u\|_{X_{\sigma}^{s_{1},s_{2},b}}
 =  \bigl\|\langle\xi\rangle^{s_{1}+\sigma} |\xi|^{-\sigma}
\langle\eta\rangle^{s_{2}}
 \left\langle\lambda\right\rangle^{b}\mathscr{F}u(\zeta)
 \bigr\|_{L_{\zeta}^{2}(\SR^{3})}.
$$
Let $X=X_{1}+X_{2}$, $\tilde{X}=\tilde{X_{1}}+\tilde{X_{2}}$,
 with
 \begin{align*}
& X_{1} =X_{0}^{s_{1},\>s_{2},\>b-b^{\prime}}, \quad
  X_{2} =X_{\sigma}^{s_{1},\>s_{2},\>b}\cap
  X_{\sigma}^{s_{1}-3\epsilon,\>s_{2},\>b+\epsilon},\nonumber\\
& \tilde{X_{1}}=X_{0}^{s_{1},\>s_{2},\>0}, \quad
   \tilde{X_{2}}=X_{\sigma}^{s_{1},\>s_{2},\>b^{\prime}}
\cap X_{\sigma}^{s_{1}-3\epsilon,\>s_{2},\>b^{\prime}+\epsilon}.
 \end{align*}
Here, $b=\frac{1}{2}+\frac{\epsilon}{2},b^{\prime}
=-\frac{1}{2}+\epsilon,$
 $\sigma=\frac{1}{2}+\epsilon$.
The norm in $X$ is given by
\begin{align*}
\|u\|_{X}={\rm inf}\left\{\|v\|_{X_{1}}+\|w\|_{X_{2}}\mid
 u=v+w,v\in X_{1},w\in X_{2}\right\}
\end{align*}
and the norm in $\tilde{X}$ is given by
\begin{align*}
\|u\|_{\tilde{X}}={\rm inf}\left\{\|v\|_{\tilde{X}_{1}}+
\|w\|_{\tilde{X}_{2}}
\mid u=v+w,v\in \tilde{X_{1}},w\in \tilde{X_{2}}\right\}.
\end{align*}
It is easy to see that (c.f. \cite{Hadac2008})
\begin{equation}\label{embedding}
  X\hookrightarrow X_{0}^{s_{1},s_{2},b}.
\end{equation}
The space $ X_{T}$ denotes the restriction
 of $X$ onto the finite time interval $[0,T]$ and
is equipped with the norm
 \begin{equation*}
    \|u\|_{X_{T}} =\inf \left\{\|u\|_{X}
    :g\in X, u(t)=g(t)
 \>\> {\rm for} \>  t\in [0,T]\right\}.
 \end{equation*}
 It is easy to see that $X_{T}\hookrightarrow
 C([-T,T];H^{s_{1},s_{2}}(\R^{2}))$.

To prove that the Cauchy problem is locally well-posed in
$H^{-\frac{1}{2},\>0}(\R^{2})$, we use the $U^{p}$ and $V^{p}$ spaces which were introduced in
\cite{KTCPAM, KTIMRN,W,Hadac2009}.

We consider functions taking values in $L^{2}=L^{2}(\R^{d},\R)$, but in the general part of this section one may replace $L^{2}$ by an arbitrary Hilbert space. Let $ \mathcal{Z}$ be the sets of finite partitions
$-\infty<t_{0}<t_{1}<\cdot\cdot\cdot<t_{K}\leq \infty$.
\begin{Definition}
Let $1\leq p<\infty$. For $\{t_{k}\}_{k=0}^{K}\in \mathcal{Z}$ and $\{\phi_{k}\}_{k=0}^{K}\subset L^{2}$
with $\sum\limits_{k=0}^{K-1}\|\phi_{k}\|_{L^{2}}^{p}=1$ and $\phi_{0}=0$, we call the function $a:\R\longrightarrow L^{2}$ given by
\begin{eqnarray*}
a=\sum\limits_{k=1}^{K}\chi_{[t_{k-1},t_{k})}\phi_{k-1}
\end{eqnarray*}
a $U^{p}$ atom. Furthermore, we define the atomic space
\begin{eqnarray*}
U^{p}=\left\{u=\sum\limits_{j=1}^{\infty}\lambda_{j}a_{j}:a_{j}\quad U^{p}-atom,\lambda_{j}\in \mathbb{C}\>s.t.\sum\limits_{j=1}^{\infty}|\lambda_{j}|<\infty\right\}
\end{eqnarray*}
endowed with the norm
\begin{eqnarray*}
\|u\|_{U^{p}}:=\inf\left\{\sum\limits_{j=1}^{\infty}|\lambda_{j}|:u=\sum\limits_{j=1}^{\infty}\lambda_{j}a_{j}\in \mathbb{C}\>s.t.\sum\limits_{j=1}^{\infty}|\lambda_{j}|<\infty\right\}.
\end{eqnarray*}
\end{Definition}
Now we present two useful statements about $U^{p}$.

\begin{Proposition}
{\rm (i)}: $\|\cdot\|_{U^{p}}$ is a norm. The space $U^{p}$ is complete and hence a Banach
space.
{\rm (ii)} The embeddings $U^{p}\subset U^{q}\subset L^{\infty}(\R,L^{2})$ are continuous.
\end{Proposition}
\setcounter{Definition}{2}
\begin{Definition}
{\rm (i)} We define $V^{p}$ as the normed space of all functions $v:\R\longrightarrow L^{2}$
such that $\lim\limits_{t\longrightarrow\pm\infty}v(t)$ exists and for which the norm
\begin{eqnarray*}
\|v\|_{V^{p}}:=\sup\limits_{\{t_{k}\}_{k=0}^{K}\in \mathcal{Z}}\left(\sum\limits_{k=1}^{K}\|v(t_{k})-v(t_{k-1})\|_{L^{2}}^{p}\right)^{\frac{1}{p}}
\end{eqnarray*}
is finite. We use the convention that $v(-\infty)=\lim\limits_{t\longrightarrow-\infty}v(t)$ and $v(\infty)=0.$
{\rm (ii)} We denote the closed subspace of all right-continuous functions $v:\R\longrightarrow L^{2}$ such that
$\lim\limits_{t\longrightarrow -\infty}v(t)=0$ by $V_{rc}^{p}.$
\end{Definition}
\setcounter{Proposition}{3}
\begin{Proposition}
Let $1\leq p< q<\infty.$
{\rm (i)}  The embedding $U^{p}\subset V_{rc}^{p}$ is continuous.
{\rm (ii)}The embedding $V^{p}\subset V^{q}$ is continuous.
{\rm (iii)}  The embedding $V_{rc}^{p}\subset U^{q}$ is continuous and
\begin{eqnarray*}
\|v\|_{U^{q}}\leq c_{p,q}\|v\|_{V^{p}}.
\end{eqnarray*}
\end{Proposition}

\begin{Proposition}
For $u\in U^{p}$ and $v\in V^{p^{\prime}}$, where $\frac{1}{p}+\frac{1}{p^{\prime}}=1$ and
a partition $\rm{t}=\{t_{k}\}_{k=0}^{K}\in \mathcal{Z}$, we define
\begin{eqnarray*}
B_{{\rm t}}(u,v):=\sum\limits \left\langle u(t_{k-1}),v(t_{k})-v(t_{k-1})\right\rangle,
\end{eqnarray*}
where
$\left\langle \cdot,\cdot \right\rangle$ denotes the inner product of $L^{2}$. Notice that
$v(t_{K})=0$ since $t_{K}=0$ for all partitions $\{t_{k}\}_{k=0}^{K}\in \mathcal{Z}$. There is a unique
number $B(u,v)$ with the property that for all $\epsilon>0$ there exists ${\rm t}\in \mathcal{Z}$ such that
for every ${\rm t^{\prime}}\subset {\rm t}$ it holds
\begin{eqnarray*}
\left|B_{{\rm t^{\prime}}}(u,v)-B(u,v)\right|<\epsilon
\end{eqnarray*}
and the associated bilinear form
\begin{eqnarray*}
B:U^{p}\times V^{p^{\prime}}:(u,v)\longrightarrow B(u,v)
\end{eqnarray*}
satisfies
\begin{eqnarray*}
\left|B(u,v)\right|\leq \|u\|_{U^{p}}\|v\|_{V^{p^{\prime}}}.
\end{eqnarray*}
\end{Proposition}
\begin{Proposition}
Let $1<p<\infty$. We have
\begin{eqnarray*}
(U^{p})^{\star}=V^{p^{\prime}}
\end{eqnarray*}
in the sense that
\begin{eqnarray*}
T:V^{p^{\prime}}\longrightarrow (U^{p})^{\star},T(v):=B(\cdot, v)
\end{eqnarray*}
is an isometric isomorphism.

\end{Proposition}
\setcounter{Corollary}{6}
\begin{Corollary}
For $1<p<\infty,u\in U^{p}$ and for $v\in V^{p}$ the following estimates hold true
\begin{eqnarray*}
\|u\|_{U^{p}}=\sup \limits_{v\in V^{p^{\prime}},\|v\|_{V^{p^{\prime}}}=1}|B(u,v)|
\end{eqnarray*}
and
\begin{eqnarray*}
\|v\|_{V^{p}}=\sup \limits_{u \>U^{p^{\prime}}-atom}|B(u,v)|.
\end{eqnarray*}
\end{Corollary}
\setcounter{Proposition}{7}
\begin{Proposition}
Let $1<p<\infty$. If the distributional derivative of $u$ is in $L^{1}$ and $v\in V^{p}$. Then,
\begin{eqnarray*}
B(u,v)=-\int_{-\infty}^{\infty}\left\langle u^{\prime}(t),v(t) \right\rangle dt.
\end{eqnarray*}

\end{Proposition}

\setcounter{Definition}{8}
\begin{Definition}
Let $s\in \R, 1\leq p,q\leq \infty$. We define the semi-norms
\begin{eqnarray*}
\|u\|_{\dot{B}_{p,q}^{s}}=\left(\sum\limits_{N}N^{qs}\|Q_{N}u\|_{L^{p}(\SR,L^{2})}^{q}\right)^{\frac{1}{p}}(q<\infty),
\|u\|_{\dot{B}_{p,\infty}^{s}}=\sup \limits_{N}N^{s}\|Q_{N}u\|_{L^{p}(\SR,L^{2})}
\end{eqnarray*}
for all $u\in \mathcal{S}^{\prime}(\R, L^{2})$ for which these numbers are finite.
\end{Definition}
\setcounter{Proposition}{9}
\begin{Proposition}
Let $1< p<\infty$. For $v\in V^{p},$ the estimate
\begin{eqnarray*}
\|v\|_{\dot{B}_{p,\infty}^{\frac{1}{p}}}\leq C\|v\|_{V^{p}}
\end{eqnarray*}
holds true.
Moreover, for any  $u\in \mathcal{S}^{\prime}(\R, L^{2})$ such that the semi-norm $\|u\|_{\dot{B}_{p,1}^{\frac{1}{p}}}$ is finite there exists $u(\pm \infty)\in L^{2}.$ Then, $u-u(-\infty)\in U^{p}$ and the estimate
\begin{eqnarray*}
\|u-u(-\infty)\|_{U^{p}}\leq C\|u\|_{\dot{B}_{p,1}^{\frac{1}{p}}}
\end{eqnarray*}
holds true.
\end{Proposition}
Now, we focus on the spatial dimension $d=2(i.e. L^{2}=L^{2}(\R^{2},\mathcal {C})$ and consider
$S=\partial_{x}^{3}+\partial_{x}^{-1}\partial_{y}^{2}-\partial_{x}^{-1}$.
We define the associated unitary operator $e^{-tS}:L^{2}\longrightarrow L^{2}$ to be the Fourier multiplier
\begin{eqnarray*}
(\mathscr{F}_{xy}e^{-t S}u_{0})(\xi,\eta)={\rm exp } \left(it(\xi^{3}-\frac{\eta^{2}+1}{\xi})\right)\mathscr{F}_{x}u_{0}(\xi,\eta).
\end{eqnarray*}
\setcounter{Definition}{10}
\begin{Definition}

{\rm (i)} $U_{S}^{p}=e^{-\cdot S}U^{p}$ with norm $\|u\|_{U_{S}^{p}}=\|e^{\cdot S}u\|_{U^{p}},$
{\rm (ii)} $V_{S}^{p}=e^{-\cdot S}V^{p}$ with norm $\|u\|_{V_{S}^{p}}=\|e^{\cdot S}u\|_{V^{p}},$
and similarly the closed subspace $U_{c,S}^{p},V_{rc, S}^{p},V_{-,S}^{p}$ and $V_{-,rc, S}^{p}.$
We define the smooth projections
\begin{eqnarray*}
\mathscr{F}P_{N}u(\xi,\eta):=\Psi_{N}(\xi)\mathscr{F}u(\xi,\eta,\tau), \mathscr{F}Q_{M}^{S}u(\xi,\eta,\tau)=
\Psi_{M}(\tau-\xi^{3}+\frac{\eta^{2}+1}{\xi})\mathscr{F}u(\xi,\eta,\tau),
\end{eqnarray*}
as well as $\mathscr{F}P_{0}u(\xi,\eta,\tau)=\chi(2\xi)\mathscr{F}u(\xi,\eta,\tau),Q_{\geq M}^{S}:=\sum\limits_{N\geq M}Q_{N}^{S}$ and $Q_{<M}^{S}=I-Q_{\geq M}^{S}.$ Obviously, $Q_{M}^{S}=e^{\cdot S}Q_{M}e^{-\cdot S}.$
\end{Definition}
\setcounter{Definition}{11}
\begin{Definition}
Let $s,b\in \R$ and $1\leq q\leq \infty$. We define the semi-norms
\begin{eqnarray*}
\|u\|_{\dot{X}^{s,b,q}}:=\left(\sum\limits_{N}N^{2s}\left\|e^{\cdot S}P_{N}u\right\|_{\dot{B}_{2,q}^{b}}^{2}\right)^{\frac{1}{2}}
\end{eqnarray*}
for all $u\in \mathcal{S}^{\prime}(\R,l^{2})$ for which these numbers are finite.Moreover,
\begin{eqnarray*}
\dot{X}^{0,\frac{1}{2},1}\subset U_{S}^{2}\subset V_{-,S}^{2}\subset \dot{X}^{0,\frac{1}{2},\infty}
\end{eqnarray*}
are continuous. We may identity $u\in \mathcal{S}^{\prime}(\R, L^{2})$ with a subset of $\mathcal {S}^{\prime}(\R^{3})$ and
\begin{eqnarray*}
\dot{X}^{s,b,q}=\left(\sum\limits_{N}N^{2s}\left(\sum\limits_{M}M^{bq}\|P_{N}Q_{M}^{S}u\|_{L^{2}(\SR^{3})}^{q}
\right)^{\frac{2}{q}}\right)^{\frac{1}{2}}
\end{eqnarray*}
with the obvious modification in the case $q=\infty.$
\end{Definition}
\setcounter{Definition}{12}
\begin{Definition}
{\rm (i)} Define $\dot{Y}^{s}$ as the closure of all $u\in C(\R; H^{1,1}(\R^{2}))\cap V_{-,rc,S}^{2}$ such that
\begin{eqnarray*}
\|u\|_{\dot{Y}^{s}}:=\left(\sum\limits_{N}N^{2s}\|P_{N}u\|_{V_{S}^{2}}^{2}\right)^{\frac{1}{2}}<\infty,
\end{eqnarray*}
in the space $C(\R;\dot{H}^{s,0}(\R^{2})$ with respect to the $\|\cdot\|_{\dot{Y}^{s}}-norm.$
{\rm (ii)} Define $\dot{Y}^{s}$ as the closure of all $u\in C(\R; H^{1,1}(\R^{2}))\cap U_{S}^{2}$ such that
\begin{eqnarray*}
\|u\|_{\dot{Z}^{s}}:=\left(\sum\limits_{N}N^{2s}\|P_{N}u\|_{U_{S}^{2}}^{2}\right)^{\frac{1}{2}}<\infty,
\end{eqnarray*}
in the space $C(\R;\dot{H}^{s,0}(\R^{2})$ with respect to the $\|\cdot\|_{\dot{Z}^{s}}-norm.$
{\rm (iii)} Define $X$ as the closure of all $u\in C(\R; H^{1,1}(\R^{2}))\cap U_{S}^{2}$ such that
\begin{eqnarray*}
\|u\|_{X}:=\|u\|_{\dot{Z}^{0}}+\|u\|_{\dot{X}^{0,1,1}}<\infty,
\end{eqnarray*}
in the space $C(\R;L^{2}(\R^{2})$ with respect to the $X-norm.$ Define $Z^{s}:=\dot{Z}^{s}+X$ with norm
\begin{eqnarray*}
\|u\|_{Z^{s}}={\rm inf}\left\{\|u_{1}\|_{\dot{Z}^{s}}+\|u_{2}\|_{X}|u=u_{1}+u_{2}\right\}.
\end{eqnarray*}
\end{Definition}

\setcounter{Definition}{13}
\begin{Definition}
Let $E$ be a Banach space of continuous functions $f:\R\longrightarrow H,$
for some Hilbert space $H$. We also
consider the corresponding restriction space to the interval $I\subset \R$ by
\begin{eqnarray*}
E(I)=\left\{u\in C(I,H)|\exists \tilde{u}\in E:\tilde{u}=u, t\in I\right\}
\end{eqnarray*}
endowed with the norm $\|u\|_{E(I)}={\rm inf}\left\{\|\tilde{u}\|_{E}|\tilde{u}:
\tilde{u}=u,t\in I\right\}$.
\end{Definition}

\setcounter{Proposition}{14}
\begin{Proposition}

{\rm (i)} Let $T>0$ and $u\in \dot{Y}^{s}([0,T])$, $u(0)=0$. Then,
for every $\epsilon>0$ there exists $0\leq T^{\prime}\leq T$ such
that $\|u\|_{\dot{Y}^{s}([0,T])}<\epsilon.$

{\rm (ii)} Let $T>0$ and $u\in \dot{Z}^{s}([0,T])$, $u(0)=0$. Then,
 for every $\epsilon>0$ there exists $0\leq T^{\prime}\leq T$ such that
  $\|u\|_{\dot{Z}^{s}([0,T])}<\epsilon.$

\end{Proposition}

The main results of this paper are in the following three theorems.

\begin{Theorem}\label{Thm1}
\mbox{\rm (Bilinear estimate)}\

\noindent Let
 $s_{1}>-\frac{1}{2},\>s_{2}\geq0$, $b=\frac{1}{2}+\frac{\epsilon}{2},
 b^{\prime}=-\frac{1}{2}+\epsilon$, and $\sigma=\frac{1}{2}+\epsilon$.
 Then, we have
$$
     \|\partial_{x}(u_{1}u_{2})\|_{\tilde{X}}
\leq C\|u_{1}\|_{X}  \> \|u_{2}\|_{X}.
$$
\end{Theorem}

\begin{Theorem}\label{Thm2}\mbox{\rm(Local well-posedness in $H^{s_{1},s_{2}}$ with $s_{1}>-\frac{1}{2},\>s_{2}\geq0.$)}\

  \noindent Let  $s_{1}>-\frac{1}{2},\>s_{2}\geq0.$  For $R>0$,
   there exists $T=T(R)>0$ such that
 for each $u_{0}\in B_{R} =\left\{u_{0}\in\! H^{s_{1},s_{2}}(\R^{2})
 \mid \|u_{0}\|_{H^{s_{1},s_{2}}(\SR^{2})}<R\right\}$,
there exists a unique solution
 to  (\ref{1.04})(\ref{1.02}) in $X_{T}.$ Moreover, the mapping
 $F_{R}:u_{0}\rightarrow u$ is analytic from $B_{R}$ to $X_{T}$.
\end{Theorem}

\begin{Theorem}\label{Thm3}\mbox{\rm (Ill-posedness in $H^{s_{1},\>0}(\R^{2})$ with $s_{1}<-\frac{1}{2}$)}\

 \noindent (\ref{1.04})-(\ref{1.02}) are ill-posed in the space
$H^{s_{1},\>0}(\R^{2})$ with $s_{1}<-\frac{1}{2}$ in
the sense that there is no $T>0$ for which
the solution map $u_{0} \longrightarrow u$
 is $C^{3}$
from $H^{s_{1},0}(\R^{2})$ to $H^{s_{1},0}(\R^{2})$ at zero.
\end{Theorem}

\begin{Theorem}\label{Thm4}\mbox{\rm (Well-posedness in $H^{-\frac{1}{2},\>0}(\R^{2})$)}\

\noindent
 For $\lambda>0$ and  all $R\geq \lambda$ and $u_{0}\in B_{\lambda,R}$ there exists a solution
 \begin{eqnarray*}
 u\in Z^{-\frac{1}{2}}([0,T])\subset C([0,T]; H^{-\frac{1}{2},0}(\R^{2}))
 \end{eqnarray*}
 for $T=\lambda^{6}R^{-6}$ of (\ref{1.04})-(\ref{1.02}) on $(0,T)$. If a solution $v\in Z^{-\frac{1}{2}}([0,T])$
 on $(0,T)$ satisfies $v(0)=u(0),$ then $v=u|_{[0,T]}.$ Moreover, the flow map
 \begin{eqnarray*}
 u_{0}(\in B_{\lambda,R})\longrightarrow u(\in Z^{-\frac{1}{2}}([0,T]))
 \end{eqnarray*}
 is analytic.
\end{Theorem}

\noindent
{\bf Remark 1.} Chen et al. \cite{CLZ} established the
local well-posedness with $s_{1}>-\frac{3}{10}$ and $s_{2}\geq 0$.
Theorem 2 considerably improves this result.
Combining Theorem 1.1 with the Banach fixed point theorem,
we can obtain Theorem 1.2. From Theorem 1.2 and  Theorem 1.3
we see that $s=-\frac{1}{2}$ is the critical exponent for well-posedness of
(\ref{1.04}) (\ref{1.02}) in $H^{s,0}(\R^{2})$.

The rest of this paper is arranged as follows. In Section 2, we present some
preliminaries. In Section 3, we derive some crucial bilinear estimates,
in particular, we prove Theorem 1.1. We prove Theorem 1.2
in Section 4 and then Theorem 1.3 in Section 5.

\bigskip

\bigskip

\setcounter{section}{2}

\noindent{\large\bf 2. Preliminaries }

\setcounter{equation}{0}

\setcounter{Theorem}{0}

\setcounter{Lemma}{0}

In this section, we  present some preliminaries, Lemmas 2.1--2.13,
  which play a significant role in establishing
  in Lemmas 3.1--3.12. Lemma 2.2 together with Lemmas 3.1--3.12
    yields Theorem 1.1.The new ingredients are  as follows. Firstly,
we split the frequency space into the regular region and the singular
region. More precisely, we define $\Omega_{0}=\Omega_{r}\cup
\Omega_{s},$ where $\Omega_{r}= \Omega_{0}\cap \Omega_{1}$ and
$\Omega_{s}= \Omega_{0}\cap \Omega_{2}$. Secondly, for $\Omega_{rA}
=\Omega_{r}\cap \Omega_{A},$ we establish Lemma 2.4; for $\Omega_{rB}
=\Omega_{r}\cap \Omega_{B},$ we establish Lemma 2.5.
Thirdly, for the singular region $\Omega_{s}$, we establish
Lemmas 2.9, 2.11, 2.13.

In this section, we always denote
\begin{align*}
  &\xi=\xi_1+\xi_2,\quad
  \eta=\eta_1+\eta_2,\quad
  \tau=\tau_1+\tau_2, \\
  &v=3\xi\xi_{1}\xi_{2},  \quad
  w=3(\xi_{1}^{2}-\xi_{1}\xi_{2}+\xi_{2}^{2}), \quad
  f(v):=v+\frac{w}{v}+\lambda.
 \end{align*}

\begin{Lemma}\label{Lemma2.1}
\mbox{\rm \cite[Lemma 2.7]{CLZ}}
Let $b> {1}/{2}$, $N=\left( {2}/{3}\right)^{1/4}$.
Then for $u\in L^{2}(\R^{3})$,
 \begin{align*}
 &\|P_{N}\mathscr{F}^{-1}
   \left(\langle \lambda \rangle ^{-\frac{11b}{12}}
    \mathscr{F}u\right)\|_{L_{xyt}^{4}}
     \leq C\|u\|_{L^{2}},\\
 &\|P^{N}
   \mathscr{F}^{-1}\left(\langle \lambda \rangle ^{-b}
    \mathscr{F}u\right)\|_{L_{xyt}^{4}}
     \leq C\|u\|_{L^{2}}.
     \end{align*}
\end{Lemma}

\begin{Lemma}\label{Lemma2.2}
 Let $b> {1}/{2}$. Then
 \begin{equation}
 \|u\|_{L_{xyt}^{4}}\leq C\|u\|_{X_{0}^{0,0,b}},
 \quad \forall\, u\in X_{0}^{0,0,b}.\label{2.01}
 \end{equation}
 \end{Lemma}
\noindent {\bf Proof.} From Lemma 2.1, we have
$$
\left\|P_{N}u\right\|_{L_{xyt}^{4}}
\leq C
\|u\|_{X_{0}^{0,\>0,\frac{11}{12}b}},\quad
\left\|P^{N}u\right\|_{L_{xyt}^{4}}
\leq C
\|u\|_{X_{0}^{0,\>0,\>b}}.
$$
This  proves Lemma 2.2.
\begin{Lemma}\label{Lemma2.3}
Let $T\in (0,1)$, $s_{1},s_{2} \in \R$ and
$- {1}/{2}<b^{\prime}\leq0\leq b\leq b^{\prime}+1$.
Then
\begin{align}
&\left\|\psi(t)W(t)\phi\right\|_{X_{0}^{s_{1},s_{2},b}}
\leq C\|\phi\|_{H^{s_{1},\>s_{2}}},\label{2.02}\\
&\left\|\psi\Big(\frac{t}{T}\Big)
  \int_{0}^{t}W(t-\tau)h(\tau)d\tau\right\|_{X_{0}^{s_{1},\,s_{2},\,b}}
\leq C T^{1+b^{\prime}-b}
  \|h\|_{X_{0}^{s_{1},\,s_{2},\,b^{\prime}}} .
\label{2.03}
\end{align}
\end{Lemma}

For the proof of Lemma 2.3, we refer the readers to
\cite{G,Bourgain93,KPV1993},
 and  Lemmas 1.7 and 1.9 of  \cite{Grunrock}.

\begin{Lemma}\label{Lemma2.4}
Let  $\Omega_{rA}=\Omega_r\bigcap A$.
Define $P_{\frac{1}{3}}(u_{1},u_{2})(\zeta)$ via
\begin{align*}
  \mathscr{F}P_{\frac{1}{3}}(u_{1},u_{2})(\zeta)
= \int_{\SR^{3}}\chi_{_{\Omega_{rA}}}
  \mathscr{F}u_{1}(\zeta_j)
  \mathscr{F}u_{2}(\zeta_2)\, d\zeta_{1} .
\end{align*}
Then, for $b>\frac{1}{2}$, we have
\begin{align}
      \big\|P_{\frac{1}{3}}(u_{1},u_{2})\big\|_{L_{txy}^{2}}
\leq C\big\||D_{x}|^{\frac{1}{2}}u_{1}\big\|_{X_{0}^{0,0,b}}\>
      \big\||D_{x}|^{-\frac{1}{2}}u_{2}\big\|_{X_{0}^{0,0,b}}.
      \label{2.04}
\end{align}
\end{Lemma}

\begin{Lemma}\label{Lemma2.5}
Let $\Omega_{rB}=\Omega_r\bigcap B$,
$P_{\frac{1}{3}}(u_{1},u_{2}) $ be as in Lemma 2.4 with
$\Omega_{rA}$ replaced by $\Omega_{rB}$. Then for $b>\frac{1}{2}$,
we have
$$
\left\|P_{\frac{1}{3}}(u_{1},u_{2})\right\|_{L_{txy}^{2}}
\leq C\left\|u_{1}\right\|_{X_{0}^{0,0,b}}
\left\||D_{x}|^{-1}u_{2}\right\|_{X_{0}^{0,0,b}}.
$$
\end{Lemma}

Lemmas 2.4 and 2.5 can be proved similarly to
\cite[Theorem 3.3]{Hadac2008}.

\begin{Lemma}\label{Lemma2.6}
Let $\Omega_{1A}=\Omega_1\bigcap A$. Then for $b>\frac{1}{2}$, we have
\begin{align}
 \left|\int_{\SR^6}\chi_{_{\Omega_{1A}}}
    \frac{|\xi_{1}|^{-\frac{1}{2}}
          |\xi_{2}|^{ \frac{1}{2}}
         }
         {
          \langle\lambda_{1}\rangle^{b}\>
          \langle\lambda_{2}\rangle^{b}
         }
         F(\zeta)
         F_{1}(\zeta_1)
         F_{2}(\zeta_2)
       d\zeta_1d\zeta \right|
  \leq C\|F\|_{L_{ \zeta}^{2}}\>
      \|F_{1}\|_{L_{ \zeta}^{2}}\>
      \|F_{2}\|_{L_{ \zeta}^{2}}.
      \label{2.05}
 \end{align}
 \end{Lemma}

Combining Lemma 2.2 with Lemma 2.4, and by using a proof similar to
Proposition 3.5 of \cite{Hadac2008}, we can derive Lemma 2.6.

Now we are going to establish the estimates in the singular region.
First we have the following lemma.

\begin{Lemma}\label{Lemma2.7} In the singular region $\Omega_r$,
when $\frac{2\sqrt{w}}{3}\leq |\lambda|\leq 4w^{\frac{4}{7}}$,
 we have
\begin{align}
\int_{\frac{\sqrt{15}}{2}|\xi|}^{2|\xi|}\langle f(v)\rangle^{-\frac{1}{2}}dv
+\int_{-2|\xi|}^{\frac{-\sqrt{15}}{2}|\xi|}
\langle f(v)\rangle^{-\frac{1}{2}}dv
\leq C\langle\lambda\rangle^{\frac{13}{16}}.
\label{2.06}
\end{align}

\end{Lemma}
\noindent{\bf Proof.}  For the  estimation of  the first integral,
we consider $\lambda\ge 0$ and $\lambda<0 $ respectively.

 When $\lambda\ge 0,$ we have
\begin{align}
\int_{\frac{\sqrt{15}}{2}|\xi|}^{2|\xi|}\langle f(v)\rangle^{-\frac{1}{2}}dv
\leq  C\int_{\frac{\sqrt{15}}{2}|\xi|}^{2|\xi|}
\langle v\rangle^{-\frac{1}{2}}dv
\leq C\langle\xi\rangle^{\frac{1}{2}}.\label{2.07}
\end{align}
When $\lambda< 0,$ we consider $|\lambda|\leq 2\sqrt{w}$ and
$|\lambda|>2\sqrt{w} $  respectively.

For the case $|\lambda|\leq 2\sqrt{w},$ it is obviously that
\begin{align}\nonumber
\int_{\frac{\sqrt{15}}{2}|\xi|}^{2|\xi|}\langle f(v)\rangle^{-\frac{1}{2}}dv
& =
\bigg(\!\int_{\frac{\sqrt{15}}{2}|\xi|}^{\sqrt{w}-|\lambda|^{\frac{3}{4}}}
+\int_{\sqrt{w}-|\lambda|^{\frac{3}{4}}}^{\sqrt{w}
+|\lambda|^{\frac{3}{4}}}
+\int_{\sqrt{w}+|\lambda|^{\frac{3}{4}}}^{2|\xi|}\!\bigg)
\langle f(v)\rangle^{-\frac{1}{2}}dv\\
&\definition
I_{1}+I_{2}+I_{3},\label{2.08}
\end{align}

It is obviously,
$
I_{2}\leq C\langle\lambda \rangle ^{\frac{3}{4}}.
$

For $\dps I_{1}$,
noting that $f(v)$ is decreasing in $ (0,\sqrt{w}]$,
when $v\in [\frac{\sqrt{15}}{2}|\xi|,\sqrt{w}-|\lambda|^{\frac{3}{4}}]$,
\begin{align}
f(v)&\geq f(\sqrt{w}-|\lambda|^{\frac{3}{4}})
=2\sqrt{w}+\lambda-|\lambda|^{\frac{3}{4}}
  +\frac{w}{\sqrt{w}-|\lambda|^{\frac{3}{4}}}-\sqrt{w}\nonumber\\
&\geq -\lambda^{\frac{3}{4}}+\frac{w}{\sqrt{w}
-|\lambda|^{\frac{3}{4}}}-\sqrt{w}=\frac{|\lambda|^{\frac{3}{2}}}
{\sqrt{w}-|\lambda|^{\frac{3}{4}}}.\label{2.09}
\end{align}
Since
$ |\lambda|\geq \frac{2\sqrt{w}}{3}
$,
from the above, we have
$f(v)\geq C|\lambda|^{\frac{1}{2}}.$
Thus, we have
\begin{align}
I_{1}\definition\!\!
\int_{\frac{\sqrt{15}}{2}|\xi|}^{\sqrt{w}-|\lambda|^{\frac{3}{4}}}\!\!
  \langle f(v)\rangle^{-\frac{1}{2}} dv
\leq C\!\int_{\frac{\sqrt{15}}{2}|\xi|}^{\sqrt{w}-|\lambda|^{\frac{3}{4}}}\!\!
 \langle\lambda\rangle^{-\frac{1}{4}} dv \leq C\sqrt{w}
\langle\lambda \rangle ^{-\frac{1}{4}}\leq C\langle
\lambda \rangle ^{\frac{3}{4}}.\label{2.010}
\end{align}

Similar to $I_{1}$, we have
$
I_{3}\leq C\langle \lambda \rangle ^{\frac{3}{4}}.
$
Therefore, when $|\lambda|\leq 2\sqrt{w},$  we have
\begin{align}
\int_{\frac{\sqrt{15}}{2}|\xi|}^{2|\xi|}\langle f(v)\rangle^{-\frac{1}{2}}dv
=I_{1}+I_{2}+I_{3}
\leq C\langle\lambda\rangle^{\frac{3}{4}}.\label{2.011}
\end{align}

Now  we consider the case $|\lambda|>2\sqrt{w}$. Let $v_{2}
=\frac{-\lambda+\sqrt{\lambda^{2}-4w}}{2}$ be the larger root of
$f(v)=v+\frac{A}{v}+\lambda=0$. Without loss of generality, we may assume
$v_{2}-|\lambda|^{\frac{3}{4}}>\sqrt{w}+|\lambda|^{\frac{3}{4}}.$

Note that $\frac{\sqrt{15}}{2}|\xi|\geq \sqrt{w},$
 we have
\begin{align}
   \int_{\frac{\sqrt{15}}{2}|\xi|}^{2|\xi|}
   \langle f(v)\rangle^{-\frac{1}{2}}dv
&\leq
   \bigg(\int_{\sqrt{w}}^{\sqrt{A}+|\lambda|^{\frac{3}{4}}}\!
  + \int_{\sqrt{w}+|\lambda|^{\frac{3}{4}}}^{v_{2}-|\lambda|^{\frac{3}{4}}}
  + \int_{v_{2}-|\lambda|^{\frac{3}{4}}}^{v_{2}+|\lambda|^{\frac{3}{4}}}
  + \int_{v_{2}+|\lambda|^{\frac{3}{4}}}^{2|\xi|}\bigg)
    \frac{dv}{\langle f(v)\rangle^{\frac{1}{2}}}\nonumber\\
&\definition
   J_{1}+J_{2}+J_{3}+J_{4}.\label{2.012}
\end{align}

 It is obvious that
\begin{align}
J_{1}\leq C\langle \lambda\rangle^{\frac{3}{4}},\quad
J_{3}\leq C|\lambda|^{\frac{3}{4}}.\label{2.013}
\end{align}

To control $J_{2}$, we consider $2\sqrt{w}\leq |\lambda|\leq
2\sqrt{w}+|\lambda|^{\frac{1}{2}}$ and
$|\lambda|>2\sqrt{w}+|\lambda|^{\frac{1}{2}}$
respectively.

 When  $2\sqrt{w}\leq |\lambda|\leq 2\sqrt{w}+|\lambda|^{\frac{1}{2}}$,
 we have
\begin{align}
&\quad\, \left|(v_{2}-|\lambda|^{\frac 34})
     -(\sqrt{w}+|\lambda|^{\frac 34})\right| \nonumber\\
&\leq  \left|\sqrt{w}-\frac{|\lambda|}{2}\right|
     + \left|\frac{|\lambda|}{2}-v_{2}\right|
     + 2|\lambda|^{\frac 34}
 \leq |\lambda|^{\frac{1}{2}}
     +\sqrt{\lambda^{2}-4w}\nonumber\\
&\leq |\lambda|^{\frac{1}{2}}+|\lambda|^{\frac{1}{4}}
      \sqrt{|\lambda|+2\sqrt{w}}
 \leq |\lambda|^{\frac{1}{2}}
       +2|\lambda|^{\frac{3}{4}}
 \leq C|\lambda|^{\frac{3}{4}}.
 \label{2.014}
\end{align}
thus we have
$
J_{2}\leq C|\lambda|^{\frac{3}{4}}.
$

When  $|\lambda|> 2\sqrt{w}+|\lambda|^{\frac{1}{2}}$,
  $f(v)$ is negative and increasing on $[\sqrt{w},+\infty)$ and since
$|\lambda|\leq 4w^{\frac{4}{7}}$, for
$v\in \Big[\sqrt{w}+|\lambda|^{\frac 34}, v_2-|\lambda|^{\frac 34}\Big]$,
we have
\begin{align}
|f(v)|
&\geq \big|f(v_{2}-|\lambda|^{\frac{3}{4}})\big|=|\lambda|^{\frac{3}{4}}
\frac{v_{2}(\sqrt{\lambda ^{2}-4w}-|\lambda|^{\frac{3}{4}})}{v_{2}(|\lambda|-|\lambda|^{\frac{3}{4}})-w}
\nonumber\\&
\geq\frac{2\sqrt{w}v_{2}|\lambda|}{[v_{2}(|\lambda|
-|\lambda|^{\frac{3}{4}})-w][\sqrt{|\lambda|+2\sqrt{w}}
+|\lambda|^{\frac{1}{2}}]}\nonumber\\
&\geq \frac{2\sqrt{w}v_{2}|\lambda|}{v_{2}|\lambda|^{\frac{3}{2}}}
\geq 2w^{\frac{1}{2}}|\lambda|^{-\frac{1}{2}}\geq |\lambda|^{\frac{3}{8}}.
\label{2.015}
\end{align}
Thus we have
\begin{align}
J_{2}\leq C|\lambda|^{\frac{13}{16}}.\label{2.016}
\end{align}
Combining (\ref{2.015}) with (\ref{2.016}) we see that for both cases,
$
J_{2}\leq C|\lambda|^{\frac{13}{16}}.
$

Similar to $J_{2}$, we have
$
J_{4}\leq C|\lambda|^{\frac{13}{16}}.
$
Therefore,
we obtain the estimate for the first integral
\begin{align}
\int_{\frac{\sqrt{15}}{2}|\xi|}^{|2|\xi|}\langle f(v)
\rangle^{-\frac{1}{2}}dv
\leq C|\lambda|^{\frac{13}{16}}.\label{2.017}
\end{align}

For the second integral in \eqref{2.06}, by  the symmetry on $f(v)$
we obtain the desired estimate.

\noindent{\bf Remark 4.} We note that in the singular region, i.e. when
$\left|1-\frac{1}{3\xi_{1}^{2}\xi_{2}^{2}}\right|< \frac{1}{4}$,
it holds $\frac{\sqrt{15}}{2}|\xi|\leq|v|\leq 2|\xi|.$

\begin{Lemma}\label{Lemma2.8}
For any $a\in \R$ and $\epsilon>0$, the following inequalities hold.
\begin{align}
&\int_{\SR}\frac{dt}{\langle t\rangle^{1+\epsilon}
\langle t-a\rangle^{1+\epsilon}}
\leq C\langle a\rangle^{-1-\epsilon},\quad
\int_{\SR}\frac{dt}
{\langle t\rangle^{1+\epsilon}|t-a|^{\frac{1}{2}}}
\leq C\langle a\rangle^{-\frac{1}{2}},\label{2.018}
\end{align}
\end{Lemma}

For the proof of this lemma, we refer to \cite[Lemma 2.8]{KPV1996}.
\begin{Lemma}\label{Lemma2.9}
For the singular region $\Omega_s$, when
$ \frac{2}{3}\sqrt{w}\leq|\lambda|\leq 4w^{\frac{4}{7}}$, we have
\begin{align}
\int_{\SR^{3}}\chi_{_{\Omega_s}}
  \frac{1}{\langle\lambda_{1}
   \rangle^{1+2\epsilon}\langle\lambda_{2}\rangle^{1+2\epsilon}}\,d\zeta_{1}
\leq C|\xi|^{-\frac{5}{2}}\>
    \langle\lambda\rangle^{\frac{13}{16}}.
\label{2.019}
\end{align}
\end{Lemma}

\noindent {\bf Proof.} By using Lemma \ref{Lemma2.8}, we have
\begin{align}
\int_{\SR^{3}}\chi_{_{\Omega_s}}
    \frac{1}
         {\langle\lambda_{1}\rangle^{1+2\epsilon}
          \langle\lambda_{2}\rangle^{1+2\epsilon}
         }\,d\zeta_{1}
\leq C\int_{\SR^{2}}\chi_{_{\Omega_s}}
       \frac{1}
       {\langle\lambda_{1}+\lambda_{2}\rangle^{1+2\epsilon}}\,
       d\xi_{1}d\eta_{1}.
       \label{2.020}
\end{align}
Now change the variables and let
$v=3\xi\xi_{1}\xi_{2},
 u=\lambda_{1}+\lambda_{2}=\tau -\phi(\xi_1,\eta_1)-\phi(\xi_2,\eta_2).$
 Then, we have
\begin{align}
\left|\frac{\partial(v,u)}{\partial(\xi_{1},\eta_{1})}\right|
=6\left|\xi_{1}^{2}-\xi_{2}^{2}\right|\left|\frac{\eta_{1}}
{\xi_{1}}-\frac{\eta_{2}}{\xi_{2}}\right|\sim |\xi|^{\frac{5}{2}}
\left|f(v)-u\right|^{\frac{1}{2}}.\label{2.021}
\end{align}
Combining (\ref{2.020}) with (\ref{2.021}), by using  Lemmas \ref{Lemma2.7}
and \ref{Lemma2.8},  we have
\begin{align*}
&\quad \,  \int_{\SR^{3}}\chi_{_{\Omega_s}}
   \frac{1}
   { \langle\lambda_{1}\rangle^{1+2\epsilon}\>
     \langle\lambda_{2}\rangle^{1+2\epsilon}}
   d\zeta_{1}\\
&\leq C\int_{\SR^{2}}\chi_{_{\Omega_s}}
     \frac{1}{|\xi|^{\frac{5}{2}}
     \left|f(v)-u\right|^{\frac{1}{2}}\langle u\rangle^{1+2\epsilon}}\,dvdu
\leq C \int_{\SR}\chi_{_{\Omega_s}}
      \frac{1}{|\xi|^{\frac{5}{2}}
       \langle f(v)\rangle^{\frac{1}{2}}}\, dv\nonumber\\
&\leq C|\xi|^{-\frac{5}{2}}
    \left[\int_{\frac{\sqrt{15}}{2}|\xi|}^{2|\xi|}
    \langle f(v)\rangle^{-\frac{1}{2}}dv
  + \int_{-2|\xi|}^{\frac{-\sqrt{15}}{2}|\xi|}
   \langle f(v)\rangle^{-\frac{1}{2}}dv{}\right]\leq C|\xi|^{-\frac{5}{2}} \langle\lambda\rangle^{\frac{13}{16}}.
\label{2.022}
\end{align*}
This completes the proof of Lemma  2.9.

\begin{Lemma}\label{Lemma2.10}
For the singular region $\Omega_s$,
when $|\lambda|\geq 4w^{\frac{4}{7}}$,
we have
\begin{equation}
     \int_{\frac{\sqrt{15}}{2}|\xi|}^{2|\xi|}
     \langle f(v)\rangle^{-\frac{1}{2}}dv
     +\int_{-2|\xi|}^{-\frac{\sqrt{15}}{2}|\xi|}
     \langle f(v)\rangle^{-\frac{1}{2}}dv
\leq C|\xi|\langle\lambda\rangle^{-\frac{1}{2}}.
      \label{2.023}
\end{equation}
\end{Lemma}
\noindent{\bf Proof.} By Remark 4 and from the assumption,
 we have
\begin{align}
\frac{\sqrt{15}}{2}|\xi|\leq |v|\leq 2|\xi|, \quad
\left|f(v)\right|\geq \frac{|\lambda|}{4}.\label{2.024}
\end{align}
Thus
$$
     \int_{\frac{\sqrt{15}}{2}|\xi|}^{2|\xi|}
     \langle f(v)\rangle^{-\frac{1}{2}}dv
     +\int_{-2|\xi|}^{-\frac{\sqrt{15}}{2}|\xi|}
     \langle f(v)\rangle^{-\frac{1}{2}}dv
\leq C|\xi|\>\langle\lambda\rangle^{-\frac{1}{2}}.
$$
This completes the proof of Lemma  2.10.

\begin{Lemma}\label{Lemma2.11}
For the singular region $\Omega_s$,
when $|\lambda|\geq 4w^{\frac{4}{7}}$,
we have
\begin{align}
\int_{\SR^{3}}\chi_{_{\Omega_s}} \frac{1}
     {\langle\lambda_{1}\rangle^{1+2\epsilon}
      \langle\lambda_{2}\rangle^{1+2\epsilon}}\,
      d \zeta_{1}
\leq C|\xi|^{-\frac{3}{2}}\>
     \langle\lambda\rangle^{-\frac{1}{2}} .
     \label{2.025}
\end{align}
\end{Lemma}
\noindent {\bf Proof.} Combining Lemma 2.10 with
a proof similar to Lemma 2.9,
 we can prove Lemma 2.11.

This completes the proof of Lemma  2.11.

\begin{Lemma}\label{Lemma2.12}
For the singular region $\Omega_s$,
when $|\lambda|< 2\sqrt{w}$,  we have
\begin{align}
   \int_{\frac{\sqrt{15}}{2}|\xi|}^{2|\xi|}
   \langle f(v)\rangle^{-\frac{1}{2}}dv
  +\int_{-2|\xi|}^{-\frac{\sqrt{15}}{2}|\xi|}
   \langle f(v)\rangle^{-\frac{1}{2}}dv
\leq C\langle\xi\rangle^{\frac{3}{4}}.\label{2.026}
\end{align}
\end{Lemma}

\noindent{\bf Proof.} We prove the estimate for the first integral.
We consider $\lambda\ge0$ and $\lambda<0,$ respectively.

  When $\lambda\ge 0,$ we have
\begin{align}
     \int_{\frac{\sqrt{15}}{2}|\xi|}^{2|\xi|}
     \langle f(v)\rangle^{-\frac{1}{2}}dv
\leq C\langle\xi\rangle^{\frac{1}{2}}.
     \label{2.027}
\end{align}

When $\lambda<0,$ since
$|\lambda|\leq 2\sqrt{w}$ and $|\xi|<\sqrt{w},$
we have
\begin{align}
    \int_{\frac{\sqrt{15}}{2}|\xi|}^{2|\xi|}
    \langle f(v)\rangle^{-\frac{1}{2}}dv
& \le \bigg(
    \int_{0}^{\sqrt{w}-|\xi|^{\frac{3}{4}}}
  + \int_{\sqrt{w}-|\xi|^{\frac{3}{4}}}^{\sqrt{w}+|\xi|^{\frac{3}{4}}}
  + \int_{\sqrt{w}+|\xi|^{\frac{3}{4}}}^{2\sqrt{w}} \bigg)
    \langle f(v)\rangle^{-\frac{1}{2}}dv\nonumber\\
&\definition I_{1}+I_{2}+I_{3}.\label{2.028}
\end{align}

It is obvious that
$
I_{2}\leq C\langle \xi \rangle ^{\frac{3}{4}}.
$

Since $f(v)$ is decreasing on $[0,\sqrt{w}]$,
for $v\in \big[0, \sqrt{w}-|\xi|^{\frac{3}{4}}\big]$ we have
\begin{align}
f(v)
&\geq f(\sqrt{w}-|\xi|^{\frac{3}{4}})
=2\sqrt{w}+\lambda-|\xi|^{\frac{3}{4}}
+\frac{w}{\sqrt{w}-|\xi|^{\frac{3}{4}}}
-\sqrt{w}\nonumber\\
&\geq -|\xi|^{\frac{3}{4}}+\frac{w}{\sqrt{w}-|\xi|^{\frac{3}{4}}}
-\sqrt{w}=\frac{|\xi|^{\frac{3}{2}}}{\sqrt{w}-|\xi|^{\frac{3}{4}}}
\geq C|\xi|^{\frac{1}{2}}.
\label{2.029}
\end{align}
Thus we have
\begin{align}
I_{1}\leq C\int_{0}^{\sqrt{w}-|\xi|^{\frac{3}{4}}}|\xi|^{-\frac{1}{4}}dv
\leq C\langle \xi \rangle ^{\frac{3}{4}}.\label{2.030}
\end{align}

Since $f(v)$ is increasing on $[\sqrt{w},+\infty)$, while
\begin{align}
&\quad\, f(\sqrt{w}+|\xi|^{\frac{3}{4}})
=2\sqrt{w}+\lambda+|\xi|^{\frac{3}{4}}
+\frac{w}{\sqrt{w}+|\xi|^{\frac{3}{4}}}
-\sqrt{w}\nonumber\\
&\geq |\xi|^{\frac{3}{4}}
      +\frac{w}{\sqrt{w}
      +|\xi|^{\frac{3}{4}}}-\sqrt{w}
=\frac{|\xi|^{\frac{3}{2}}}
      {\sqrt{w} +|\xi|^{\frac{3}{4}}}
 \geq C|\xi|^{\frac{1}{2}}>0,\label{2.031}
\end{align}
Thus, for $v\in \big[ \sqrt{w}+|\xi|^{\frac{3}{4}} ,  2\sqrt{w}\big]$,
 $f(v)\geq f(\sqrt{w}+|\xi|^{\frac{3}{4}})\ge C|\xi|^{\frac{1}{2}} $.
 Therefore  we have
$
I_{3}\leq C\langle \xi \rangle ^{\frac{3}{4}}.
$

For the second integral in \eqref{2.025}, by  the symmetry on $f(v)$
 we obtain the desired estimate.

The proof of Lemma 2.12 is completed.

\begin{Lemma}\label{Lemma2.13}
For the singular region $\Omega_s$,
when $|\lambda|< 2\sqrt{w}$,  we have
$$
\int_{\SR^{3}}\chi_{_{\Omega_s}}
 \frac{1}
{\langle\lambda_{1}\rangle^{1+2\epsilon}
\langle\lambda_{2}\rangle^{1+2\epsilon}}d \zeta_{1}
\leq C\langle\xi\rangle^{-\frac{7}{4}}.
$$
\end{Lemma}

\noindent {\bf Proof.} Combining Lemma 2.12 with a proof
similar to Lemma 2.9,
we see that Lemma 2.13 is valid.

\begin{Lemma}\label{Lemma2.14}
We have
\begin{eqnarray}
&&\left\|Q_{M}^{S}u\right\|_{L^{2}(\SR^{3})}\leq CM^{-\frac{1}{2}}\|u\|_{V_{S}^{2}},\label{2.031}\\
&&\left\|Q_{\geq M}^{S}u\right\|_{L^{2}(\SR^{3})}\leq CM^{-\frac{1}{2}}\|u\|_{V_{S}^{2}},\label{2.032}\\
&&\left\|Q_{<M}^{S}u\right\|_{V_{S}^{p}}\leq C\|u\|_{V_{S}^{p}},\left\|Q_{\geq M}^{S}u\right\|_{V_{S}^{p}}\leq C\|u\|_{V_{S}^{p}}\label{2.033}\\
&&\left\|Q_{<M}^{S}u\right\|_{U_{S}^{p}}\leq C\|u\|_{U_{S}^{p}},\left\|Q_{\geq M}^{S}u\right\|_{U_{S}^{p}}\leq C\|u\|_{U_{S}^{p}}\label{2.034}.
\end{eqnarray}
\end{Lemma}

 Lemma 2.14 can be proved similarly to Corollary 2.18 of \cite{Hadac2008}.

\begin{Lemma}\label{Lemma2.15}
Let $T_{0}: L^{2}\times \cdot\cdot \cdot \times L^{2}\longrightarrow L_{loc}^{1}(\R^{2}; C)$ be a $n$-linear operator.

\noindent (i) Assume that for some $1\leq p,q\leq \infty$
\begin{eqnarray}
\left\|T_{0}(e^{\cdot S}\phi_{1},\cdot\cdot\cdot,e^{\cdot S}\phi_{n})\right\|_{L_{t}^{p}(\SR;L^{q}_{xy}(\SR^{2})}\leq C\prod\limits_{i=1}^{n}\|\phi_{i}\|_{L^{2}}.
\end{eqnarray}
Then, there exists $T:U_{S}^{p}\times \cdot\cdot\cdot \times U_{S}^{p}\longrightarrow L_{t}^{p}(\R;L^{q}_{xy}(\R^{2})$ satisfying
\begin{eqnarray}
\left\|T(u_{1},\cdot\cdot \cdot u_{n})\right\|_{L_{t}^{p}(\SR;L^{q}_{xy}(\SR^{2})}\leq C\prod\limits_{i=1}^{n}\|\phi_{i}\|_{U_{S}^{p}}
\end{eqnarray}
such that $T(u_{1},\cdot \cdot \cdot ,u_{n})(t)(x,y)=T_{0}(u_{1},\cdot \cdot \cdot ,u_{n})(x,y) a.e.$

\noindent (ii)
 Assume that for some $1\leq p,q\leq \infty$
\begin{eqnarray}
\left\|T_{0}(e^{\cdot S}\phi_{1},\cdot\cdot\cdot,e^{\cdot S}\phi_{n})\right\|_{L_{x}^{q}(\SR;L^{p}_{ty}(\SR^{2})}\leq C\prod\limits_{i=1}^{n}\|\phi_{i}\|_{L^{2}}.
\end{eqnarray}
Then, for $r=: {\rm min}(p,q),$ there exists $T:U_{S}^{r}\times \cdot\cdot\cdot \times U_{S}^{r}\longrightarrow L_{x}^{q}(\R;L^{p}_{ty}(\R^{2})$ satisfying
\begin{eqnarray}
\left\|T(u_{1},\cdot\cdot \cdot u_{n})\right\|_{L_{x}^{q}(\SR;L^{p}_{ty}(\SR^{2})}\leq C\prod\limits_{i=1}^{n}\|\phi_{i}\|_{U_{S}^{r}}
\end{eqnarray}
such that $T(u_{1},\cdot \cdot \cdot ,u_{n})(t)(x,y)=T_{0}(u_{1},\cdot \cdot \cdot ,u_{n})(x,y) a.e.$

\end{Lemma}

Lemma 2.15 can be proved similarly to Proposition  2.19 of \cite{Hadac2008}.
\begin{Lemma}\label{Lemma2.16}
Let $q>1$, $E$ be a Banach space and $T:U_{S}^{p}\longrightarrow E$ be a bounded, linear operator with
$\left\|Tu\right\|_{E}\leq C_{q}\|u\|_{U_{S}^{q}}$  for all $u\in U_{S}^{Q}$. In addition, assume that
for some $1\leq p <q$ there exists $C_{p}\in (0,C_{q})$ such that the estimate $\|Tu\|_{E}\leq C_{p}\|u\|_{U_{S}^{p}}$.
Then, $T$ satisfies the estimate
\begin{eqnarray}
\|Tu\|_{E}\leq \frac{4C_{p}}{\alpha_{p,q}}\left(In \frac{C_{q}}{C_{p}}+2\alpha_{p,q}+1\right)\|u\|_{V_{S}^{p}},u\in V_{-,\>rc,\>S}^{p},
\end{eqnarray}
where
$\alpha _{p,q}=(1-\frac{p}{q})In2.$
\end{Lemma}

Lemma 2.16 can be proved similarly to Proposition  2.20 of \cite{Hadac2008}.

\begin{Lemma}\label{Lemma2.17}
We have
\begin{eqnarray}
&&\left\|u\right\|_{L^{4}(\SR^{3})}\leq C\|u\|_{U_{S}^{4}},\label{2.040}\\
&&\left\|u\right\|_{L^{4}(\SR^{3})}\leq C\|u\|_{V_{-,S}^{p}}(1\leq p<4),\label{2.041}\\
&&\left\|\partial_{x}u\right\|_{L_{x}^{\infty}(\SR;L_{ty}^{2}(\SR^{2}))}\leq C\|u\|_{U_{S}^{2}},\label{2.042}.
\end{eqnarray}
\end{Lemma}

 Lemma 2.17 can be proved similarly to (24)-(26) of Corollary 2.21 of \cite{Hadac2008}.

 \begin{Lemma}\label{Lemma2.18}
For $b>\frac{1}{2},$ we have
\begin{align}
 \left|\int_{\SR^6}\chi_{_{\Omega_{rA}}}
    \frac{\mathscr{F}u\prod\limits_{j=1}^{2}\mathscr{F}P_{N_{j}}u_{j}}
         {
          \langle\lambda_{1}\rangle^{b}\>
          \langle\lambda_{2}\rangle^{b}
         }d\zeta_1d\zeta \right|
  \leq C\left(\frac{N_{1}}{N_{2}}\right)^{\frac{1}{2}}\|\mathscr{F}u\|_{L_{ \zeta}^{2}}\>
      \prod\limits_{j=1}^{2}\|\mathscr{F}P_{N_{j}}u_{j}\|_{L_{ \zeta}^{2}}\>
      .
      \label{2.043}
 \end{align}
 In particular, let $\mathscr{F}P_{N_{j}}u_{j}$ possess the same support  as in (\ref{2.043}),  we have
 \begin{eqnarray}
&&\left\|P_{N_{1}}uP_{N_{2}}u_{2}\right\|_{L^{2}(\SR^{3})}\leq C\left(\frac{N_{1}}{N_{2}}\right)^{\frac{1}{2}}\|u_{1}\|_{U_{S}^{2}}\|u_{2}\|_{U_{S}^{2}}\label{2.044}.
\end{eqnarray}
and
\begin{eqnarray}
\left\|P_{N_{1}}uP_{N_{2}}u_{2}\right\|_{L^{2}(\SR^{3})}\leq C\left(\frac{N_{1}}{N_{2}}\right)^{\frac{1}{2}}\left(In \frac{N_{2}}{N_{1}}+1\right)^{2}\|u_{1}\|_{V_{S}^{2}}\|u_{2}\|_{V_{S}^{2}}\label{2.045}.
\end{eqnarray}
\end{Lemma}

Combining Lemma 2.5 with (27)-(28) of Corollary 2.21 of \cite{Hadac2008}, we have that Lemma 2.18 is valid with the aid of  Lemma 2.16.

 \begin{Lemma}\label{Lemma2.19}
For $b>\frac{1}{2}$, we have
\begin{align}
 \left|\int_{\SR^6}\chi_{_{\Omega_{rB}}}
    \frac{\mathscr{F}u\prod\limits_{j=1}^{2}\mathscr{F}P_{N_{j}}u_{j}}
         {
          \langle\lambda_{1}\rangle^{b}\>
          \langle\lambda_{2}\rangle^{b}
         }d\zeta_1d\zeta \right|
  \leq C \left(\frac{N_{1}}{N_{2}}\right)^{\frac{1}{2}}\|\mathscr{F}u\|_{L_{ \zeta}^{2}}\>
      \prod\limits_{j=1}^{2}\|\mathscr{F}P_{N_{j}}u_{j}\|_{L_{ \zeta}^{2}}\>
      .
      \label{2.046}
 \end{align}
 In particular, let $P_{N_{j}}u_{j}$ possess the same support as in (\ref{2.046}),  we have
 \begin{eqnarray}
&&\left\|P_{N_{1}}uP_{N_{2}}u_{2}\right\|_{L^{2}(\SR^{3})}\leq C\left(\frac{N_{1}}{N_{2}}\right)^{\frac{1}{2}}\|u_{1}\|_{U_{S}^{2}}\|u_{2}\|_{U_{S}^{2}}\label{2.047}.
\end{eqnarray}
and
\begin{eqnarray}
\left\|P_{N_{1}}uP_{N_{2}}u_{2}\right\|_{L^{2}(\SR^{3})}\leq C\left(\frac{N_{1}}{N_{2}}\right)^{\frac{1}{2}}\left(In \frac{N_{1}}{N_{2}}+1\right)^{2}\|u_{1}\|_{V_{S}^{2}}\|u_{2}\|_{V_{S}^{2}}\label{2.048}.
\end{eqnarray}
\end{Lemma}

Combining Lemma 2.6 with (27)-(28) of Corollary 2.21 of \cite{Hadac2008}, we have that Lemma 2.19 is valid with the aid of Lemma 2.16.

 \begin{Lemma}\label{Lemma2.20}
For $b>\frac{1}{2}$ and $\frac{2}{3}\sqrt{w}\leq |\lambda|\leq 4|w|^{\frac{4}{7}}$, we have
\begin{align}
 \left|\int_{\SR^6}\chi_{_{\Omega_{s}}}
    \frac{\mathscr{F}u\prod\limits_{j=1}^{2}\mathscr{F}P_{N_{j}}u_{j}}
         {
          \langle\lambda_{1}\rangle^{b}\>
          \langle\lambda_{2}\rangle^{b}
         }d\zeta_1d\zeta \right|
  \leq CN_{2}^{-\frac{5}{4}}\langle \lambda \rangle^{\frac{13}{32}}\|\mathscr{F}u\|_{L_{ \zeta}^{2}}\>
      \prod\limits_{j=1}^{2}\|\mathscr{F}P_{N_{j}}u_{j}\|_{L_{ \zeta}^{2}}\>
      .
      \label{2.049}
 \end{align}
 In particular, let $P_{N_{j}}u_{j}$ possess the same support as in (\ref{2.049}),  we have
 \begin{eqnarray}
&&\left\|P_{N_{1}}uP_{N_{2}}u_{2}\right\|_{L^{2}(\SR^{3})}\leq C\left(\frac{1}{N_{2}}\right)^{\frac{11}{14}}\|u_{1}\|_{U_{S}^{2}}\|u_{2}\|_{U_{S}^{2}}\label{2.050}.
\end{eqnarray}
and
\begin{eqnarray}
\left\|P_{N_{1}}uP_{N_{2}}u_{2}\right\|_{L^{2}(\SR^{3})}\leq C\left(\frac{1}{N_{2}}\right)^{\frac{11}{14}}\left(In N_{2}+1\right)^{2}\|u_{1}\|_{V_{S}^{2}}\|u_{2}\|_{V_{S}^{2}}\label{2.051}.
\end{eqnarray}
\end{Lemma}
 \noindent
 {\bf Proof.} By using the Cauchy-Schwarz inequality, we have
\begin{align*}
&\left|\int_{\SR^{6}}
    \chi_{\Omega_s}
   \frac{\mathscr{F}u\prod\limits_{j=1}^{2}\mathscr{F}P_{N_{j}}u_{j}}{\langle\lambda_{1}\rangle^{b}\>
          \langle\lambda_{2}\rangle^{b}}
        d\zeta_1d\zeta\right|
\leq&\, C
 \int_{\SR^{3}}
    I(\zeta)
   \mathscr{F}u
    \left[\int_{\SR^{3}}
    \chi_{\Omega_s}
    \prod\limits_{j=1}^{2}\mathscr{F}P_{N_{j}}u_{j}
    d\zeta_{1} \right]^{\frac{1}{2}}
    d\zeta,
\end{align*}
where
\begin{align}
  I(\zeta)
=
  \left[\int_{\SR^{3}}
    \chi_{\Omega_s}
    \frac{1}
   {\langle\lambda_{1}\rangle^{2b}\,
    \langle\lambda_{2}\rangle^{2b}}
     d\zeta_{1} \right]^{\frac{1}{2}}.
      \label{2.052}
\end{align}
 By Lemma 2.9 and  $\frac{2}{3}\sqrt{w}\leq |\lambda|\leq 4|w|^{\frac{4}{7}}$,  we have
\begin{align}
I(\zeta)\leq C
|N_{2}|^{-\frac{5}{4}}\langle \lambda\rangle^{\frac{13}{32}}\leq CN_{2}^{-\frac{11}{14}}.\label{2.053}
\end{align}
Combining (\ref{2.053}) with the Cauchy-Schwarz inequality, we have that (\ref{2.049}) is valid.
In particular, we have that (\ref{2.050}) is  valid.
 We obtain that (\ref{2.051}) is valid with the aid of Lemma 2.16.

This completes the proof of Lemma 2.20.

 \begin{Lemma}\label{Lemma2.21}
For $b>\frac{1}{2}$ and $|\lambda|\geq 4|w|^{\frac{4}{7}}$, we have
\begin{align}
 \left|\int_{\SR^6}\chi_{_{\Omega_{s}}}
    \frac{\mathscr{F}u\prod\limits_{j=1}^{2}\mathscr{F}P_{N_{j}}u_{j}}
         {
          \langle\lambda_{1}\rangle^{b}\>
          \langle\lambda_{2}\rangle^{b}
         }d\zeta_1d\zeta \right|
  \leq CN_{2}^{-\frac{3}{4}}\langle \lambda \rangle^{-\frac{1}{4}}\|\mathscr{F}u\|_{L_{ \zeta}^{2}}\>
      \prod\limits_{j=1}^{2}\|\mathscr{F}P_{N_{j}}u_{j}\|_{L_{ \zeta}^{2}}\>
      .
      \label{2.054}
 \end{align}
 In particular, let $P_{N_{j}}u_{j}$ possess the same support as in (\ref{2.054}),  we have
 \begin{eqnarray}
&&\left\|P_{N_{1}}uP_{N_{2}}u_{2}\right\|_{L^{2}(\SR^{3})}\leq C\left(\frac{1}{N_{2}}\right)^{\frac{29}{28}}\|u_{1}\|_{U_{S}^{2}}\|u_{2}\|_{U_{S}^{2}}\label{2.055}.
\end{eqnarray}
and
\begin{eqnarray}
\left\|P_{N_{1}}uP_{N_{2}}u_{2}\right\|_{L^{2}(\SR^{3})}\leq C\left(\frac{1}{N_{2}}\right)^{\frac{29}{28}}\left(In N_{2}+1\right)^{2}\|u_{1}\|_{V_{S}^{2}}\|u_{2}\|_{V_{S}^{2}}\label{2.056}.
\end{eqnarray}
\end{Lemma}

Combining Lemma 2.11 with a proof similar to Lemma 2.20, we have that Lemma 2.21 is valid.
\bigskip

\bigskip

\noindent{\large\bf 3. Bilinear estimates}

\setcounter{equation}{0}

 \setcounter{Theorem}{0}

\setcounter{Lemma}{0}

 \setcounter{section}{3}
 In this section, we  present the proof of some bilinear
estimates which are the core of this paper. The most difficult
case is the interaction between the high frequency and the low frequency
 in establishing bilinear estimates. More precisely,  the establishment of
  Lemma
 3.2-3.7 is most difficult. We only present details of proof of
  Lemmas 3.2-3.5,
 other cases can be similarly proved. Now we present the new ingredients
  in establishing Lemmas 3.2-3.5. Firstly, we
   divide the integration domain into three parts:
\begin{equation*}
\mbox{\rm (i)\, Region $\Omega_{r }\bigcap A$;\qquad
 (ii)\, Region $\Omega_r\bigcap B$;\qquad
 (iii)\, Region $\Omega_s$} .
\end{equation*}
Secondly, for part (i), we use Lemma 2.4 to establish bilinear estimates;
for part (ii), we use Lemma 2.5 to establish bilinear estimates;
for part (iii), we use Lemmas 2.9, 2.11, 2.13 to establish bilinear estimates.

To establish the desired bilinear estimates, we define first the
following operator (c.f. \cite{Hadac2008})
 \begin{align*}
\mathscr{F}P_{c}(u_{1},u_{2})(\zeta)
=\int_{\SR^{3}}
  \chi_{_{|\xi_{1}|\leq c|\xi_{2}|}}
   \mathscr{F}u_{1}(\zeta_1)\,
    \mathscr{F}u_{2}(\zeta_2)\, d \zeta_{1} ,
\end{align*}

Let $c=1$. We decompose the operator further as follows.
\begin{align*}
\partial_{x}P_{1}(u_{1},u_{2})=Q_{00}(u_{1},u_{2})
+\sum\limits_{k=1}^{2}
 \sum\limits_{j=0}^{2}Q_{kj}(u_{1},u_{2}),
\end{align*}
where the operators $Q_{k j}$ are defined by
\begin{align*}
   \mathscr{F}Q_{k j}(u_{1},u_{2})
=  i\xi \int_{\SR^{3}}
    \chi_{_{A_{k j}}}(\zeta_1, \zeta)
    \mathscr{F} u_{1}(\zeta_{1})\>
    \mathscr{F} u_{2}(\zeta-\zeta_{1})
    d\zeta_{1},
\end{align*}
where
$$
A_{00}:=\left\{(\zeta_1,\zeta)
\in \R^{6}\,\big|\, |\xi_{1}|\leq |\xi-\xi_{1}|
\leq 6^{7}\right\},
$$
and $A_{k j}:=\Xi_{k}\cap \Lambda_{j}$ for
$1\leq k\leq 2,\,0\leq j\leq2$,
and
\begin{align*}
&\Xi_{1}=\Big\{(\zeta_1,\zeta)\in \R^{6}
  \, \Big|\, |\xi_{1}|\leq \frac{1}{3}|\xi-\xi_{1}|,|\xi-\xi_{1}|
  \geq 6^{7}\Big\},\nonumber\\
&\Xi_{2}=\Big\{(\zeta_1,\zeta)\in \R^{6}
    \, \Big|\,\frac{1}{3}|\xi-\xi_{1}|\leq |\xi_{1}|
    \leq |\xi-\xi_{1}|,|\xi-\xi_{1}|\geq 6^{7}\Big\},\nonumber\\
&\Lambda_{0}=\big\{(\zeta_1,\zeta)\in \R^{6}
     \, \big|\,|\lambda|=|\lambda_{\rm max}|\big\},\nonumber\\
&\Lambda_{j}=\big\{(\zeta_1,\zeta)\in \R^{6}
     \, \big|\,|\lambda_{j}|=|\lambda_{\rm max}|\big\}\ \ (j=1,2).
\end{align*}
We always write
\begin{align*}
  \xi=\xi_1+\xi_2,\quad
  \eta=\eta_1+\eta_2,\quad
  \tau=\tau_1+\tau_2.
\end{align*}

\begin{Lemma}\label{Lemma3.1}
Let $s_{1}\geq-\frac{1}{2}+8\epsilon,s_{2}\geq0$,
$b=\frac{1}{2}+\frac{\epsilon}{2}$.
Then, for any  $u_{j}\in
 X_{0}^{s_{1},s_{2},b}$ $(j=1,2)$, we have
\begin{align}
\|Q_{00}(u_{1},u_{2})\|_{X_{0}^{s_{1},s_{2},0}}\leq C
\|u_{1}\|_{X_{0}^{s_{1},s_{2},b}}\>\|u_{2}\|_{X_{0}^{s_{1},s_{2},b}}.
\label{3.01}
\end{align}
\end{Lemma}
\noindent{\bf Proof.}  By duality,
it suffices to  prove the following inequality
\begin{align}
 \left|\int_{\SR^{3}}\bar{u}
    Q_{00}(u_{1},u_{2})dxdydt\right|
\le C\|u\|_{X_{0}^{-s_{1},\,-s_{2},\,0}}
     \|u_{1}\|_{X_{0}^{s_{1},\, s_{2},\,b}}
     \|u_{2}\|_{X_{0}^{s_{1},\,s_{2},\,b}}\label{3.02}
\end{align}
hold for any $u\in X_{0}^{-s_{1},-s_{2},0}.$
We write
\begin{align}
&  G_{0}(\zeta)\definition
    \langle\xi\rangle^{-s_{1}}
    \langle\eta\rangle^{-s_{2}}
    \mathscr{F}u(\zeta), \label{G0}\\
&  F_{j}(\zeta_j)\definition
   \langle\xi_{j}\rangle^{s_1}
   \langle\eta_{j}\rangle^{s_2}
   \langle \lambda_{j}\rangle^{b}
   \mathscr{F}u_{j}(\zeta_j)\>(j=1,2),\label{Fj} \\
&  K_{0}(\zeta_1,\zeta)\definition
   \frac{
        |\xi|\langle\xi\rangle^{s_{1}}
        \langle \eta\rangle ^{s_{2}}
        }
        {
   \langle \xi_{1}\rangle^{ s_{1}}
   \langle \eta_{1}\rangle^{ s_{2}}
   \langle \lambda_{1}\rangle^{ b}
   \langle \xi_{2}\rangle^{ s_{1}}
   \langle \eta_{2}\rangle^{ s_{2}}
   \langle \lambda_{2}\rangle^{ b}
        },
\end{align}
then (\ref{3.02}) is equivalent  to
\begin{align}
\left|\int_{\SR^{6}}
    \chi_{_{A_{00}}}
    K_{0}(\zeta_{1},\zeta)
    G_{0}(\zeta)
    F_{1}(\zeta_{1})
    F_{2}(\zeta_{2})
    d\zeta_1d\zeta\right|
\leq C \|G_{0}\|_{L_{\zeta}^{2}}\>
       \|F_{1}\|_{L_{\zeta}^{2}}\>
       \|F_{2}\|_{L_{\zeta}^{2}}.
       \label{3.03}
\end{align}
By the definition of $Q_{00}$, or $A_{00}$, one has
$|\xi_{1}|\le |\xi_{2}|\le  1 $, and
$\langle \eta\rangle ^{s_{2}}
\leq
\langle \eta_{1}\rangle ^{s_{2}}
\langle \eta_{2}\rangle ^{s_{2}},$ and thus
\begin{align}
     K_{0}(\zeta_1,\zeta)
\leq \frac{C}{
     \langle \lambda_{1}\rangle^{b}\,
     \langle \lambda_{2}\rangle^{b}}.\label{3.04}
\end{align}
Putting (\ref{3.04}) into \eqref{3.03}, applying the
H\"older inequality and the Plancherel identity as well
as Lemma 2.2, we have
\begin{align*}
&\quad \, \left|\int_{\SR^{6}}
   \chi_{_{A_{00}}}
   K_{0}(\zeta_1,\zeta)
   G_{0}(\zeta)
   F_{1}(\zeta_1)
   F_{2}(\zeta_2)
   d\zeta_1d\zeta\right|
   \nonumber\\
&\leq C\left|\int_{\SR^{6}}\frac{1}
    { \langle \lambda_{1}\rangle^{b} \>
     \langle \lambda_{2}\rangle^{b}}
     G_{0}(\zeta)\,
     F_{1}(\zeta_1)\,
     F_{2}(\zeta_2)
     d\zeta_1d\zeta\right|
  \nonumber\\
&\leq C\|G_{0}\|_{L_{ \zeta }^{2}}
  \left\|\int_{\SR^{3}}\frac{1}
     {\langle \lambda_{1}\rangle^{b} \>
      \langle \lambda_{2}\rangle^{b}}
     F_{1}(\zeta_1)\>
     F_{2}(\zeta_2)
     d\zeta_{1} \right\|_{L_{ \zeta }^{2}}\nonumber\\
&\leq C\|G_{0}\|_{L_ {\zeta}^{2}}\>
     \left\|\mathscr{F}^{-1}
     \left( \langle \lambda_{1}\rangle ^{-b}F_{1}
     \right)\right\|_{L_{xyt}^{4}}\>
     \left\| \mathscr{F}^{-1}
     \left( \langle \lambda_{2}\rangle ^{-b}F_{2}
     \right)\right\|_{L_{xyt}^{4}}\nonumber\\
&\leq C\|G_{0}\|_{L_ {\zeta}^{2}}\>
       \|F_{1}\|_{L_ {\zeta}^{2}}\>
       \|F_{2}\|_{L_ {\zeta}^{2}}.
\end{align*}
The proof of Lemma 3.1 is completed.

\begin{Lemma}\label{Lemma3.2}
Let $s_{1}\geq-\frac{1}{2}+8\epsilon,s_{2}\geq0$,
$b=\frac{1}{2}+\frac{\epsilon}{2}$,
$b^{\prime}=-\frac{1}{2}+\epsilon,$
$\sigma=\frac{1}{2}+\epsilon$.
Then, for any  $u_{j}\in X_{0}^{s_ 1 ,\,s_{2},\,b}$ $(j=1,2)$,
we have
\begin{align}
 \|Q_{10}(u_{1}u_{2})\|_{X_{\sigma}
                          ^{s_{1}-3\epsilon,\,s_{2},\,
                          b^{\prime}+\epsilon}}
\leq C \|u_{1}\|_{X_{0}^{s_{1},s_{2},b}}\>
       \|u_{2}\|_{X_{0}^{s_{1},s_{2},b}}.
     \label{3.05}
\end{align}
\end{Lemma}
\noindent{\bf Proof.}
By duality, it suffices to  prove
\begin{align}
 \left|\int_{\SR^{3}}
  \bar{u} Q_{10}(u_{1},u_{2})
   dxdydt\right|
\leq C\|u\|_{X_{-\sigma}
              ^{-s_{1}+3\epsilon,\, -s_{2},\,
                -b^{\prime}-\epsilon}}\>
      \|u_{1}\|_{X_{0}^{s_{1},\,s_{2},\,b}}\>
      \|u_{2}\|_{X_{0}^{s_{1},\,s_{2},\,b}}
      \label{3.06}
\end{align}
for any
$u\in X_{-\sigma}^{-s_{1}+3\epsilon,\, -s_{2},\, -b^{\prime}-\epsilon}.$
We write
\begin{align*}
   G_{1}(\zeta)
& \definition
    \langle\xi\rangle^{-s_{1}+3\epsilon-\sigma}\>
    |\xi|^{\sigma}\>
    \langle\eta\rangle^{-s_{2}}
    \langle\lambda\rangle^{-b^{\prime}-\epsilon}
    \mathscr{F}u(\zeta), \\
   K_{1}
& \definition
  K_{1}(\zeta_1,\zeta)
= \frac{ |\xi|^{1-\sigma}\>
         \langle \xi \rangle^{s_{1}-3\epsilon+\sigma}\>
         \langle \eta \rangle ^{s_{2}}
       }
       {
         \langle \lambda \rangle^{-b^{\prime}-\epsilon}\>
         \langle \xi_{1}\rangle^{s_{1}}\>
         \langle \eta_{1}\rangle^{s_{2}}\>
         \langle \lambda_{1}\rangle^{b} \>
         \langle \xi_{2}\rangle^{s_{1}}\>
         \langle \eta_{2}\rangle^{s_{2}}\>
         \langle \lambda_{2}\rangle^{b}
       },
\end{align*}
$F_{j}(\zeta_j), j=1,2,$ are the same as in \eqref{Fj}.
Then (\ref{3.06}) is equivalent to
\begin{align}
& \left|\int_{\SR^{6}}
     \chi_{_{A_{10}}}
      K_{1}(\zeta_1,\zeta)\>
      G_{1}(\zeta)\>
      F_{1}(\zeta_1)\>
      F_{2}(\zeta_2)\,
      d\zeta_1d\zeta\right|
\leq  C \|G_{1}\|_{L_{\zeta}^{2}}\>
        \|F_{1}\|_{L_{\zeta}^{2}}\>
        \|F_{2}\|_{L_{\zeta}^{2}}.\label{3.07}
\end{align}
We divide the integration domain into three parts:
\begin{equation}
\mbox{\rm (i)\, Region $\Omega_{r }\bigcap A$;\qquad
 (ii)\, Region $\Omega_r\bigcap B$;\qquad
 (iii)\, Region $\Omega_s$} .\label{3Parts}
\end{equation}
\noindent
We remark that (i) and (ii) are away from the singular region,
while (iii) is  a singular region part.

In Part (i), one has
$|\xi_{2}|\geq 6^{7}$, $|\xi_{1}|\leq \frac{1}{3}|\xi_2|$
and  $\langle \eta\rangle ^{s_{2}}\leq
\langle \eta_{1}\rangle ^{s_{2}}\,
\langle \eta_{2}\rangle ^{s_{2}}$, thus
\begin{align}
     K_{1}(\zeta_1,\zeta)
\leq \frac
     {
      C|\xi|^{\epsilon}\,
      \langle \xi_{1}\rangle^{-s_{1}}\,
      |\xi_{1}|^{-\frac{1}{2}+2\epsilon}
     }
     {
     \langle \lambda_{1}\rangle^{b}\,
     \langle \lambda_{2}\rangle^{b}}
 \leq C\frac{|\xi_{1}|^{-\frac{1}{2}}
             |\xi_{2}|^{\frac{1}{2}-5\epsilon}}
     {
     \langle \lambda_{1}\rangle^{b}\,
     \langle \lambda_{2}\rangle^{b}
     }
\leq C\frac{|\xi_{1}|^{-\frac{1}{2}}
            |\xi_{2}|^{\frac{1}{2}}}
     {
     \langle \lambda_{1}\rangle^{b}\,
     \langle \lambda_{2}\rangle^{b}}.\label{3.010}
\end{align}
Putting (\ref{3.010}) into \eqref{3.07}, applying the
Plancherel identity and the H\"older inequality as
well as Lemma 2.5, we have
\begin{align*}
&  \left|\int_{\SR^6}
    \chi_{_{A_{10}\cap A\cap \Omega_1}}
    K_{1}(\zeta_1,\zeta)\,
    G_{1}(\zeta)\,
    F_{1}(\zeta_1)\,
    F_{2}(\zeta_2)\,
    d\zeta_1d\zeta\right| \\
\leq & C \bigg|\int_{\SR^{6}}
    \chi_{_{A_{10}\cap A\cap \Omega_1}}
     \frac{|\xi_{1}|^{-\frac{1}{2}}\,
           |\xi_{2}|^{\frac{1}{2}}}
       {\langle \lambda_{1}\rangle^{b}\,
        \langle \lambda_{2}\rangle^{b}}
    G_{1}(\zeta)\,
    F_{1}(\zeta_1)\,
    F_{2}(\zeta_2)\,
      d\zeta_1d\zeta\bigg|\\
\leq& C\|G_{1}\|_{L_{\zeta}^{2}}\,
      \|F_{1}\|_{L_{\zeta}^{2}}\,
      \|F_{2}\|_{L_{\zeta}^{2}}.
\end{align*}
In Part (ii), $|\xi_{1}|\leq1$. Consequently,
\begin{align*}
     K_{1}(\zeta_1,\zeta)
\leq \frac{C|\xi_{2}|}
     {\langle \lambda_{1}\rangle^{b}\,
      \langle \lambda_{2}\rangle^{b}},
\end{align*}
combining the above inequality with Lemma 2.4, we have
\begin{align}\nonumber
&\quad\,
    \left|\int_{\SR^{6}}
    \chi_{_{A_{10}\cap B\cap \Omega_1}}
    K_{1}(\zeta_1,\zeta)\,
    G_{1}(\zeta)\,
    F_{1}(\zeta_1)\,
    F_{2}(\zeta_2)\,
    d\zeta_1d\zeta\right|\\
&\le \left|\int_{\SR^{6}}
     \chi_{_{A_{10}\cap B\cap \Omega_1}}
     \frac{|\xi_{2}|}
    {\langle \lambda_{1}\rangle^{b}\,
     \langle \lambda_{2}\rangle^{b}}\,
      G_{1}(\zeta)\,
      F_{1}(\zeta_1)\,
      F_{2}(\zeta_2)\,
      d\zeta_1d\zeta \right|\nonumber\\
&\leq C\|G_{1}\|_{L_{\zeta}^{2}}
       \left\|\int_{\SR^{3}}\chi_{_{A_{10}\cap B\cap \Omega_1}}
       \frac{|\xi_{2}|}
      {\langle \lambda_{1}\rangle^{b}\,
       \langle \lambda_{2}\rangle^{b}}\,
      F_{1}(\zeta_1)\,
      F_{2}(\zeta_2)\,d\zeta_1
       \right\|_{L_{\zeta}^{2}} \nonumber\\
&\leq C\|G_{1}\|_{L_{\zeta}^{2}}\,
       \|F_{1}\|_{L_{\zeta}^{2}}\,
       \|F_{2}\|_{L_{\zeta}^{2}}.
       \label{Part(ii)}
\end{align}
In Part (iii),
$|\xi_{1}|\sim |\xi_{2}|^{-1}\sim |\xi|^{-1}$. Thus,
\begin{align*}
      K_{1}(\zeta_1,\zeta)
\leq \frac{C|\xi|^{1-3\epsilon}}
     {
     \langle \lambda \rangle ^{-b^{\prime}-\epsilon}\,
     \langle\lambda_{1}\rangle^{b}\,
     \langle\lambda_{2}\rangle^{b}
     },
\end{align*}
consequently,
\begin{align*}
&\quad\,
   \left|\int_{\SR^6}
    \chi_{_{A_{10}\cap \Omega_s}}
    K_{1}(\zeta_1,\zeta)\,
    G_{1}(\zeta)\,
    F_{1}(\zeta_1)\,
    F_{2}(\zeta_2)\,
    d\zeta_1d\zeta\right|\\
&\leq C \int_{\SR^{3}}
    I_{1}(\zeta)
    G_{1}(\zeta)
    \left[\int_{\SR^{3}}
    \chi_{_{A_{10}\cap  \Omega_s}}
    |F_{1}(\zeta_1)|^{2}\,
    |F_{2}(\zeta_2)|^{2}\,
    d\zeta_{1} \right]^{\frac{1}{2}}
    d\zeta,
\end{align*}
where
\begin{align}
   I_{1}(\zeta)
=  \frac{|\xi|^{1-3\epsilon}}
   {\langle \lambda \rangle ^{-b^{\prime}-\epsilon}}
  \left[\int_{\SR^{3}}
    \chi_{_{A_{10}\cap  \Omega_s}}\frac{1}
   {\langle\lambda_{1}\rangle^{2b}\,
    \langle\lambda_{2}\rangle^{2b}}
     d\zeta_{1}
       \right]^{\frac{1}{2}}.\label{3.011}
\end{align}
Note that $0<\epsilon<\frac{1}{100}$,
$b^{\prime}=-\frac{1}{2}+\epsilon$,
and $|\lambda|\geq C|\xi|,$  using Lemmas 2.9, 2.11
and 2.13,   we have
\begin{align}
I_{1}(\zeta)\leq\frac{|\xi|^{1-3\epsilon}}
{\langle \lambda \rangle ^{-b^{\prime}-\epsilon}}{\rm max}\left\{
|\xi|^{-\frac{5}{4}}\langle \lambda\rangle^{\frac{13}{32}},
|\xi|^{-\frac{3}{4}}
\langle \lambda\rangle^{-\frac{1}{4}},
\langle\xi\rangle^{-\frac{7}{8}}\right\}\leq C.\label{3.012}
\end{align}
By using (\ref{3.011})-(\ref{3.012}) and the Cauchy-Schwartz inequality,
we have that
\begin{align}\nonumber
&   \left|\int_{\SR^{6}}
    \chi_{_{A_{10}\cap  \Omega_s}}
    K_{1}(\zeta_1,\zeta)\,
    G_{1}(\zeta)\,
    F_{1}(\zeta_1)\,
    F_{2}(\zeta_2)\,
    d\zeta_1d\zeta\right|\\
\leq& C\|G_{1}\|_{L_{\zeta}^{2}}\,
      \|F_{1}\|_{L_{\zeta}^{2}}\,
      \|F_{2}\|_{L_{\zeta}^{2}}.
      \label{Part(iii)}
\end{align}

The proof of Lemma 3.2 is completed.

\begin{Lemma}\label{Lemma3.3}
Let $s_{1}\geq-\frac{1}{2}+8\epsilon,s_{2}\geq0$,
$b=\frac{1}{2}+\frac{\epsilon}{2}$,
$b^{\prime}=-\frac{1}{2}+\epsilon,$
$\sigma=\frac{1}{2}+\epsilon$.
Then, for any
$
 u_{1} \in X_{\sigma}^{s_{1},s_{2},b},
 u_{2} \in X_{0}^{s_{1},s_{2},b}
$, we have
\begin{align}
      \|Q_{10}(u_{1},u_{2})\|_{X_{\sigma}^{s_{1},\,  s_{2},\, b^{\prime}}}
\leq C\|u_{1}\|_{X_{\sigma}^{s_{1},\,s_{2},\, b}}\,
      \|u_{2}\|_{X_{0}^{s_{1},\, s_{2},\, b}}.
\label{3.013}
\end{align}
\end{Lemma}
\noindent{\bf Proof.}
By duality, it suffices to prove
\begin{align}
    \bigg|\int_{\SR^{3}}\bar{u}Q_{10}(u_{1},u_{2})dxdydt
    \bigg|
\leq C\|u\|_{X_{-\sigma}^{-s_{1},\, -s_{2},\, -b^{\prime}}}
      \|u_{1}\|_{X_{\sigma}^{s_{1},\, s_{2},\, b}}\>
      \|u_{2}\|_{X_{0}^{s_{1},\,s_{2},\, b}}
\label{3.014}
\end{align}
for any $u\in X_{-\sigma}^{-s_{1},-s_{2},-b^{\prime}}.$
We write
\begin{align}
   G_{2}(\zeta)
&\definition
     \langle\xi\rangle^{-s_{1}-\sigma}
     |\xi|^{\sigma}\langle\eta\rangle^{-s_{2}}
     \langle\lambda\rangle^{-b^{\prime}}
     \mathscr{F}u(\zeta),\label{G2} \\
   F_{1}(\zeta_1)
&\definition
     \langle\xi_{1}\rangle^{s_1+\sigma}
     |\xi_1|^{-\sigma}
     \langle\eta_{1}\rangle^{s_2}
     \langle \lambda_{1}\rangle^{b}
     \mathscr{F}u_{1}(\zeta_1)   \\
   K_{2}(\zeta_1,\zeta)
&\definition\frac{|\xi|^{1-\sigma}
     \langle\xi\rangle^{s_{1}+\sigma}
     \langle \eta\rangle ^{s_{2}}}
     {\langle \lambda\rangle^{-b^{\prime}}
   \langle \xi_{1}\rangle^{s_{1}}
   \langle \eta_{1}\rangle^{s_{2}}
   \langle \lambda_{1}\rangle^{b}
   \langle \xi_{2}\rangle^{s_{1}}
   \langle \eta_{2}\rangle^{s_{2}}
   \langle \lambda_{2}\rangle^{b}},\label{K2}
\end{align}
and $F_{2}$ is the same as \eqref{Fj}. Then (\ref{3.014}) is equivalent to
\begin{align}
&\left|\int_{\SR^{6}}
     \chi_{_{A_{10}}}
    K_{2}(\zeta_1,\zeta)\,
    G_{2}(\zeta)\,
    F_{1}(\zeta_1)\,
    F_{2}(\zeta_2)\,
    d\zeta_1d\zeta\right|
\leq   C\|G_{2}\|_{L_{\zeta}^{2}}\,
       \|F_{1}\|_{L_{\zeta}^{2}}\,
       \|F_{2}\|_{L_{\zeta}^{2}}. \label{3.22}
\end{align}

As in the proof of Lemma 3.2, we divide the integration
domain of \eqref{3.22} also into three parts as \eqref{3Parts}.

In Part (i),
 $|\xi_{2}|\geq 6^{7}$, $|\xi_{1}|\leq \frac{1}{3}|\xi_2|$ and
 $\langle \eta\rangle ^{s_{2}}\leq
\langle \eta_{1}\rangle ^{s_{2}}\,
\langle \eta_{2}\rangle ^{s_{2}}$. Note that $b^{\prime}=-\frac{1}{2}+\epsilon$
and $s_{1}\geq-\frac{1}{2}+8\epsilon,$ we have
\begin{align}
     K_{2}(\zeta_1,\zeta)
\leq C \frac{|\xi|\,
       \langle\xi_{1}\rangle^{-s_{1}}}
      {\langle \lambda\rangle^{-b^{\prime}}
      \langle \lambda_{1}\rangle^{b}\,
       \langle \lambda_{2}\rangle^{b}}
\leq C \frac{|\xi|^{1+2b^{\prime}}\,
       |\xi_{1}|^{b^{\prime}}\,
       \langle \xi_{1}\rangle^{-s_{1}}
       }
       {\langle \lambda_{1}\rangle^{b}\,
       \langle \lambda_{2}\rangle^{b}}
\leq C\frac{|\xi_{1}|^{-\frac{1}{2}}\,
            |\xi_{2}|^{\frac{1}{2}}}
      {\langle \lambda_{1}\rangle^{b}\,
       \langle \lambda_{2}\rangle^{b}},\label{3.016}
\end{align}
the same bound as (\ref{3.010}). Thus this part can be obtained in the
same way as Part (i) in the proof of Lemma 3.2.

Part (ii)  can be proved similarly to \eqref{Part(ii)}  of Lemma 3.2.

In Part (iii),
$|\xi_{1}|\sim |\xi_{2}|^{-1}\sim |\xi|^{-1}$,
thus
\begin{align*}
     K_{2}(\zeta_1,\zeta)
\leq
    \frac{C|\xi|}
    {\langle \lambda \rangle ^{-b^{\prime}}
     \langle \lambda_{1}\rangle^{b}\,
     \langle \lambda_{2}\rangle^{b}}.
\end{align*}
Consequently, we have
\begin{align*}
&\left|\int_{\SR^{6}}
    \chi_{_{A_{10}\cap  \Omega_s}}
    K_{2}(\zeta_1,\zeta)\,
    G_{2}(\zeta)\,
    F_{1}(\zeta_1)\,
    F_{2}(\zeta_2)\,
    d\zeta_1d\zeta\right|\\
\leq&\, C
 \int_{\SR^{3}}
    I_{2}(\zeta)
    G_{2}(\zeta)
    \left[\int_{\SR^{3}}
    \chi_{_{A_{10}\cap  \Omega_s}}
    |F_{1}(\zeta_1)|^{2}\,
    |F_{2}(\zeta_2)|^{2}\,
    d\zeta_{1} \right]^{\frac{1}{2}}
    d\zeta,
\end{align*}
where
\begin{align}
   I_{2}(\zeta)
= \frac{|\xi|}
   {\langle \lambda \rangle ^{-b^{\prime}}}
  \left[\int_{\SR^{3}}
    \chi_{_{A_{10}\cap  \Omega_s}}
    \frac{1}
   {\langle\lambda_{1}\rangle^{2b}\,
    \langle\lambda_{2}\rangle^{2b}}
     d\zeta_{1} \right]^{\frac{1}{2}}.
      \label{3.017}
\end{align}
 Note that $0<\epsilon<\frac{1}{100}$
 and $b^{\prime}=-\frac{1}{2}+\epsilon,$
 by Lemmas 2.9, 2.11  and 2.13,  we have
\begin{align}
I_{2}(\zeta)\leq C\frac{|\xi|}{\langle \lambda \rangle
^{-b^{\prime}}}{\rm max}\left\{
|\xi|^{-\frac{5}{4}}\langle \lambda\rangle^{\frac{13}{32}},
|\xi|^{-\frac{3}{4}}
\langle \lambda\rangle^{-\frac{1}{4}},\langle\xi\rangle^{-\frac{7}{8}}
\right\}\leq C.\label{3.018}
\end{align}
Similar to \eqref{Part(iii)} we obtain the estimate for this part.

The proof of Lemma 3.3 is completed.

\begin{Lemma}\label{Lemma3.4}
Let $s_{1}\geq-\frac{1}{2}+8\epsilon,s_{2}\geq0$,
$b=\frac{1}{2}+\frac{\epsilon}{2}$,
$b^{\prime}=-\frac{1}{2}+\epsilon,$ $\sigma=\frac{1}{2}+\epsilon$.
Then, for any $u_{1}\in X_{0}^{s_{1},s_{2},b-b^{\prime}}$ and
 $u_{2}\in X_{0}^{s_{1},s_{2},b}$,
we have
\begin{align}
     \|Q_{11}(u_{1},u_{2})\|
        _{X_{0} ^{s_{1},\, s_{2},\, 0}}
\leq C \|u_{1}\|_{X_{0}^{s_{1},s_{2},\,b-b^{\prime}}}\>
       \|u_{2}\|_{X_{0}^{s_{1},\,s_{2},\, b}}.
       \label{3.019}
\end{align}
\end{Lemma}
\noindent{\bf Proof.}
By duality, it suffices to  prove
\begin{align}
     \left|\int_{\SR^{3}} \bar{u}
     Q_{11}(u_{1},u_{2})dxdydt\right|
\le  C\|u\|_{X_{0}^{-s_{1},-s_{2},0}}\,
      \|u_{1}\|_{X_{0}^{s_{1},s_{2},b-b^{\prime}}}\,
      \|u_{2}\|_{X_{0}^{s_{1},s_{2},b}}\label{3.020}
\end{align}
for any $u\in X_{0}^{-s_{1},-s_{2},\,0}.$
We write
\begin{align*}
    F_{1}(\zeta_1)
&=  \langle\xi_{1}\rangle^{s_{1}}
    \langle\eta_{1}\rangle^{s_{2}}\,
    \langle \lambda_{1}\rangle^{b-b^{\prime}}\,
    \mathscr{F}u_{1}(\zeta_1),\nonumber\\
    K_{3}(\zeta_1,\zeta)
&=  \frac{|\xi|\,
    \langle\xi\rangle^{s_{1}}\,
    \langle \eta\rangle ^{s_{2}}}
    {\langle \lambda_{1}\rangle^{b-b^{\prime}}\,
     \langle \lambda_{2}\rangle^{b}\,
     \langle \xi_{1}\rangle^{s_{1}}\,
     \langle \eta_{1}\rangle^{s_{2}}\,
     \langle \xi_{2}\rangle^{s_{1}}\,
     \langle \eta_{2}\rangle^{s_{2}}},
\end{align*}
$G_{0}(\zeta)$ and $F_{2}(\zeta_2)$ are the same as
\eqref{G0} and \eqref{Fj}, then (\ref{3.020}) is
equivalent  to
\begin{align}
&\left|\int_{\SR^{6}}
    \chi_{_{A_{11}}}
    K_{3}(\zeta_1,\zeta)\,
    G_{0}(\zeta)\,
    F_{1}(\zeta_1)\,
    F_{2}(\zeta_2)\,
    d\zeta_1d\zeta\right|
\leq   C\|G_{0}\|_{L_{\zeta}^{2}}\,
       \|F_{1}\|_{L_{\zeta}^{2}}\,
       \|F_{2}\|_{L_{\zeta}^{2}}. \label{3.28}
\end{align}

As in the proof of Lemma 3.2, we divide the
integration domain of \eqref{3.28} also into three parts as
\eqref{3Parts}.

In Part (i),
 $|\xi_{2}|\geq 6^{7}$, $|\xi_{1}|\leq \frac{1}{3}|\xi_2|$ and
 $\langle \eta\rangle ^{s_{2}}\leq
\langle \eta_{1}\rangle ^{s_{2}}\,
\langle \eta_{2}\rangle ^{s_{2}}$. Note that
$b^{\prime}=-\frac{1}{2}+\epsilon$
and $s_{1}\geq-\frac{1}{2}+8\epsilon,$ we have
\begin{align}
      K_{3}(\zeta_1,\zeta)
\leq \frac{C|\xi|
      \langle \xi_{1}\rangle^{-s_{1}}}
     {\langle \lambda_{1}\rangle^{b-b^{\prime}}\,
      \langle \lambda_{2}\rangle^{b}}
\leq C\frac{|\xi_{2}|^{2\epsilon}
      \langle\xi_{1}\rangle^{-s_{1}}\,
      |\xi_{1}|^{-\frac{1}{2}+2\epsilon}}
      {\langle \lambda_{1}\rangle^{b-b^{\prime}}\,
      \langle \lambda_{2}\rangle^{b}}
\leq C\frac{|\xi_{1}|^{-\frac{1}{2}}\,
      |\xi_{2}|^{\frac{1}{2}}}
      {\langle \lambda_{1}\rangle^{b}\,
      \langle \lambda_{2}\rangle^{b}}.\label{3.022}
\end{align}
a similar bound as (\ref{3.010}). Thus this part can
be obtained similarly to Part (i) in the proof of Lemma 3.2.

Part (ii)  can be proved similarly to \eqref{Part(ii)}  of Lemma 3.2.

In Part (iii), $|\xi_{1}|
\sim |\xi_{2}|^{-1}\sim |\xi|^{-1}$. Thus we have
\begin{align}
     K_{3}(\zeta_1,\zeta)
\leq \frac{C|\xi|}
     {\langle \lambda_{1}\rangle^{b-b^{\prime}}
      \langle \lambda_{2}\rangle^{b}}.
      \label{3.023}
\end{align}
Consequently, by the Cauchy-Schwarz inequality, we obtain
\begin{align}
&\left|\int_{\SR^{6}}
    \chi_{_{A_{11}\cap \Omega_s}}
    K_{3}(\zeta_1,\zeta)\,
    G_{0}(\zeta)\,
    F_{1}(\zeta_1)\,
    F_{2}(\zeta_2)\,
    d\zeta_1d\zeta\right|\nonumber\\
\leq \, & C
    \int_{\SR^{3}}
    I_{3}(\zeta)
    G_{0}(\zeta)
    \left[\int_{\SR^{3}}
    \chi_{_{A_{11}\cap  \Omega_s}}
    |F_{1}(\zeta_1)|^{2}\,
    |F_{2}(\zeta_2)|^{2}\,
    d\zeta_{1} \right]^{\frac{1}{2}}
    d\zeta,
    \label{3.024}
\end{align}
where
\begin{align}
   I_{3}(\zeta)
=  |\xi|
  \left[\int_{\SR^{3}}
    \chi_{_{A_{11}\cap  \Omega_s}}
    \frac{1}
   {\langle\lambda_{1}\rangle^{2(b-b')}\,
    \langle\lambda_{2}\rangle^{2b}}
     d\zeta_{1} \right]^{\frac{1}{2}}.
      \label{3.025}
\end{align}

We estimate $ I_{3}(\zeta)$ in thee cases respectively:

\centerline{
(a) \ $ 2\sqrt{w} \le |\lambda|\le 4w^{\frac{4}{7}}$; \qquad
(b) \ $ |\lambda| >  4w^{\frac{4}{7}}$; \qquad
(c) \ $|\lambda|< 2\sqrt{w}$.
}

In Case (a), by  Lemma 2.9 and noting $b^{\prime}
=-\frac{1}{2}+\epsilon$, we have
\begin{align}
      I_{3}(\zeta)
\leq C \frac{|\xi|}
      {\langle\lambda\rangle^{-b'}}
      \left[\int_{\SR^{3}}
      \chi_{_{A_{11}\cap  \Omega_s}}
      \frac{1}
      {\langle\lambda_{1}\rangle^{2b}\,
       \langle\lambda_{2}\rangle^{2b}}
      d\zeta_{1} \right]^{\frac{1}{2}}
\leq C|\xi|^{-\frac{1}{4}}\,
     \langle\lambda\rangle ^{b^{\prime}+\frac{13}{32}}
\leq C.\label{3.I3a}
\end{align}

 In Case (b),  by Lemma 2.11,
noting $0<\epsilon<\frac{1}{100}$, we have
\begin{align}
      I_{3}(\zeta)
&\leq C\frac{|\xi|}
      {\langle\lambda\rangle^{-b'}}
      \left[\int_{\SR^{3}}
      \chi_{_{A_{11}\cap  \Omega_s}}
      \frac{1}
      {\langle\lambda_{1}\rangle^{2b}\,
       \langle\lambda_{2}\rangle^{2b}}
      d\zeta_{1}\right]^{\frac{1}{2}}\nonumber\\
&\leq C|\xi|^{\frac{1}{4}}\,
      \langle \lambda\rangle ^{b^{\prime}-\frac{1}{4}}
\leq C|\xi|^{-\frac{17}{28}+\frac{8}{7}\epsilon}
\leq C.\label{3.I3b}
\end{align}

In Case (c),  by  Lemma 2.13, noting
$0<\epsilon<\frac{1}{100}$ and $b^{\prime}=
-\frac{1}{2}+\epsilon,$ we have
\begin{align}
      I_{3}(\zeta)
\leq C|\xi| ^{\frac 12  +\epsilon}
      \left[\int_{\SR^{3}}
      \chi_{_{A_{11}\cap  \Omega_s}}
      \frac{1}
      {\langle\lambda_{1}\rangle^{2b}\,
       \langle\lambda_{2}\rangle^{2b}}
      d\zeta_{1}\right]^{\frac{1}{2}}
\leq C|\xi|^{-\frac{3}{8}+\epsilon}
\leq C.\label{3.I3c}
\end{align}
From (\ref{3.I3a})--(\ref{3.I3c}) we get $ I_{3}(\zeta)
\leq C$. Putting this into \eqref{3.024}, applying
the Cauchy-Schwarz inequality, we have
\begin{align*}
&\quad\,
   \left|\int_{\SR^{6}}
   \chi_{_{A_{11}\cap  \Omega_s}}
    K_{3}(\zeta_1,\zeta)\,
    G_{0}(\zeta)\,
    F_{1}(\zeta_1)\,
    F_{2}(\zeta_2)\,
    d\zeta_1d\zeta\right|\nonumber\\
&\leq  C\int_{\SR^{3}}G_{0}(\zeta)
    \left[\int_{\SR^{3}}
    \chi_{_{A_{11}\cap  \Omega_s}}
    |F_{1}(\zeta_1)|^{2}\,
    |F_{2}(\zeta_2)|^{2}\,
    d\zeta_{1} \right]^{\frac{1}{2}}
    d\zeta\\
&\leq   C\|G_{0}\|_{L_{\zeta}^{2}}\,
       \|F_{1}\|_{L_{\zeta}^{2}}\,
       \|F_{2}\|_{L_{\zeta}^{2}}.
\end{align*}

This completes the proof of Lemma 3.4.

\begin{Lemma}\label{Lemma3.5}
Let $s_{1}\geq-\frac{1}{2}+8\epsilon,s_{2}\geq0$,
 $b=\frac{1}{2}+\frac{\epsilon}{2}$,
$b^{\prime}=-\frac{1}{2}+\epsilon,$ $\sigma=\frac{1}{2}+\epsilon$.
Then for any $u_{1}\in X_{\sigma}^{s_{1}-3\epsilon,s_{2},b+\epsilon}$
and $u_{2}\in X_{0}^{s_{1},s_{2},b}$,
we have
\begin{align}
       \|Q_{11}(u_{1},u_{2})\|_{X_{\sigma}
            ^{s_{1}-3\epsilon,s_{2},b^{\prime}+\epsilon}}
\leq C \|u_{1}\|_{X_{\sigma}^{s_{1}-3\epsilon,s_{2},b+\epsilon}}\,
       \|u_{2}\|_{X_{0}^{s_{1},s_{2},b}}.\label{3.027}
\end{align}
\end{Lemma}
\noindent{\bf Proof.}
By duality, it suffices to  prove
\begin{align}
       \left|\int_{\SR^{3}}\bar{u}
           Q_{11}(u_{1},u_{2})dxdydt\right|
\leq C \|u\|_{X_{-\sigma}
              ^{-s_{1}+3\epsilon,-s_{2},-b^{\prime}-\epsilon}}\,
       \|u_{1}\|_{X_{\sigma}^{s_{1}-3\epsilon,s_{2},b+\epsilon}}\,
       \|u_{2}\|_{X_{0}^{s_{1},s_{2},b}}
       \label{3.37}
\end{align}
for any
$u\in X_{-\sigma}^{-s_{1}+3\epsilon,\,-s_{2},\,-b^{\prime}-\epsilon}.$
We write
\begin{align}
    G_{3}(\zeta)
&=  \langle\xi\rangle^{-s_{1}-\sigma+3\epsilon}\,
    |\xi|^{\sigma}\,
    \langle\eta\rangle^{-s_{2}}\,
    \langle\lambda\rangle^{-b^{\prime}-\epsilon}\,
    \mathscr{F}u(\zeta),\label{G3}\\
    F_{1}(\zeta_1)
&=  \langle\xi_{1}\rangle^{s_{1}+\sigma-3\epsilon}\,
    |\xi_{1}|^{-\sigma}\,
    \langle\eta_{1}\rangle^{s_{2}}\,
    \langle \lambda_{1}\rangle^{b+\epsilon}\,
    \mathscr{F}u_{1}(\zeta_1),\nonumber\\
    K_{4}(\zeta_1,\zeta)
&=  \frac{\langle\xi\rangle^{s_{1}-3\epsilon+\sigma}\,
           |\xi|^{1-\sigma}\,
           \langle \eta\rangle ^{s_{2}}\,
           |\xi_{1}|^{\sigma}}
         { \langle \lambda\rangle^{-b^{\prime}-\epsilon}\,
           \langle\xi_{1}\rangle^{s_{1}-3\epsilon+\sigma}\,
           \langle \eta_{1}\rangle^{s_{2}}\,
           \langle \lambda_{1}\rangle^{b+\epsilon}\,
           \langle \xi_{2}\rangle^{s_{1}}\,
           \langle \eta_{2}\rangle^{s_{2}}\,
           \langle \lambda_{2}\rangle^{b}
         },\label{K4}
\end{align}
 $F_{2}(\zeta_2)$ is the same as in \eqref{Fj}.
Then \eqref{3.37} is equivalent to
\begin{align}
& \left|\int_{\SR^{6}}
      \chi_{_{A_{11}}}
      K_{4}(\zeta_1,\zeta)\>
      G_{3}(\zeta)\>
      F_{1}(\zeta_1)\>
      F_{2}(\zeta_2)\,
      d\zeta_1d\zeta\right|
\leq  C \|G_{3}\|_{L_{\zeta}^{2}}\>
        \|F_{1}\|_{L_{\zeta}^{2}}\>
        \|F_{2}\|_{L_{\zeta}^{2}}.\label{3.40}
\end{align}

Analogous to the proof of Lemma 3.2, we divide the integration domain
\eqref{3.40} also into three parts as \eqref{3Parts}.

In Part (i), i.e. in $\Omega_r\cap A$, one has
$|\xi_{2}|\geq 6^{7}$, $|\xi_{1}|\leq \frac{1}{3}|\xi_2|$
and  $\langle \eta\rangle ^{s_{2}}\leq
\langle \eta_{1}\rangle ^{s_{2}}$
$\langle \eta_{2}\rangle ^{s_{2}}$, thus
$$
     K_{4}(\zeta_1,\zeta)
\leq \frac{C|\xi|^{1-3\epsilon}|\xi_{1}|^{\sigma}}
     {\langle\lambda\rangle^{-b^{\prime}-\epsilon}\,
      \langle \xi_{1}\rangle^{s_{1}-3\epsilon+\sigma}\,
      \langle \lambda_{1}\rangle^{b+\epsilon}\,
      \langle \lambda_{2}\rangle^{b}}.
$$
We argue in two cases: $|\xi_{1}|\geq1$ and $|\xi_{1}|<1$.

Using the fact that
$
\langle \lambda\rangle^{b^{\prime}+\epsilon}
 \langle\lambda_{1}\rangle^{-b-\epsilon}
  \leq
   \langle\lambda\rangle^{-b}
    \langle \lambda_{1}\rangle^{b^{\prime}}
$
and noting that $b^{\prime}=-\frac{1}{2}+\epsilon$,  $s_{1}\geq
-\frac{1}{2}+8\epsilon,$
for both cases  we have
\begin{align}
     K_{4}(\zeta_1,\zeta)
\leq \frac{C}
     {\langle \lambda\rangle^{b}
      \langle \lambda_{2}\rangle^{b}}.
\label{3.030}
\end{align}
Thus this  part can be proved similarly to Lemma 3.1.

For Part (ii), one has $|\xi_1|\le 1$, and thus
\begin{align*}
     K_{4}(\zeta_1,\zeta)
\leq \frac{C|\xi_2|}
     {\langle \lambda\rangle^{b}
      \langle \lambda_{2}\rangle^{b}}.
\end{align*}
Thus this part can be proved similarly to \eqref{Part(ii)} of Lemma 3.2.

In Part (iii), one has
$
     |\xi_{1}|\sim |\xi_{2}|^{-1}
\sim |\xi|^{-1}
$.
Noting also that
$\langle \lambda\rangle^{b^{\prime}+\epsilon}
 \langle\lambda_{1}\rangle^{-b-\epsilon}
  \leq
   \langle\lambda\rangle^{-b}
    \langle \lambda_{1}\rangle^{b^{\prime}},
$
$
\sigma=\frac{1}{2}+\epsilon
$,
we  have
\begin{align}\nonumber
      K_{4}(\zeta_1,\zeta)
&\leq \frac{C|\xi|^{1-3\epsilon}\,
            |\xi_{1}|^{\sigma}
          }
          { \langle\lambda \rangle ^{-b^{\prime}-\epsilon}\,
            \langle \lambda_{1}\rangle^{b+\epsilon}\,
            \langle \lambda_{2}\rangle^{b}}\\
&\leq \frac{C|\xi|^{1-3\epsilon}\,
            |\xi_{1}|^{\sigma}
          }
          {\langle\lambda\rangle ^{-b^{\prime}}\,
           \langle \lambda_1\rangle^{b}\,
           \langle \lambda_{2}\rangle^{b}}
\leq \frac{C|\xi|^{\frac{1}{2}-2\epsilon-\sigma}}
          {\langle \lambda\rangle^{b}
           \langle \lambda_{2}\rangle^{b}}
\leq \frac{C}
          {\langle \lambda\rangle^{b}
           \langle \lambda_{2}\rangle^{b}}.
           \label{3.031}
\end{align}
Thus this  part can be proved similarly to Lemma 3.1.

This completes the proof of Lemma 3.5.

\begin{Lemma}\label{Lemma3.6}
Let $s_{1}\geq-\frac{1}{2}+8\epsilon,s_{2}\geq0$,
$b=\frac{1}{2}+\frac{\epsilon}{2}$,
$b^{\prime}=-\frac{1}{2}+\epsilon,$ $\sigma=\frac{1}{2}+\epsilon$.
Then for any
$u_{1}\in X_{\sigma}^{s_{1}-3\epsilon,s_{2},b+\epsilon}$
and $u_{2}\in X_{0}^{s_{1},s_{2},b}$,
we have
\begin{align}
     \|Q_{11}(u_{1},u_{2})\|_{X_{\sigma}^{s_{1},s_{2},b^{\prime}}}
\leq C\|u_{1}\|_{X_{\sigma}^{s_{1}-3\epsilon,s_{2},b+\epsilon}}\,
      \|u_{2}\|_{X_{0}^{s_{1},s_{2},b}}.\label{3.032}
\end{align}
\end{Lemma}
Lemma 3.6 can be proved similarly to Lemma 3.5.

\begin{Lemma}\label{Lemma3.7}
Let $s_{1}\geq-\frac{1}{2}+8\epsilon,s_{2}\geq0$,
$b=\frac{1}{2}+\frac{\epsilon}{2}$,
$b^{\prime}=-\frac{1}{2}+\epsilon,$ $\sigma=\frac{1}{2}+\epsilon$.
Then for any
$u_{j}\in X_{0}^{s_{1},s_{2},b}\,(j=1,2)$,
we have
\begin{align}
&    \|Q_{12}(u_{1},u_{2})\|_{X_{\sigma}^{s_{1}-3\epsilon,s_{2},
             b^{\prime}+\epsilon}}
\leq C\|u_{1}\|_{X_{0}^{s_{1},s_{2},b}}\,
      \|u_{2}\|_{X_{0}^{s_{1},s_{2},b}},
      \label{3.033}\\
&     \|Q_{12}(u_{1},u_{2})\|_{X_{\sigma}^{s_{1},s_{2},b^{\prime}}}
\leq C\|u_{1}\|_{X_{0}^{s_{1},s_{2},b}}\,
      \|u_{2}\|_{X_{0}^{s_{1},s_{2},b}}.
      \label{3.034}
\end{align}
\end{Lemma}
\noindent
(\ref{3.033})-(\ref{3.034})  can be proved similarly to Lemmas 3.2, 3.3,
 respectively.

\begin{Lemma}\label{Lemma3.8}
Let $s_{1}\geq-\frac{1}{2}+8\epsilon,s_{2}\geq0$,
$b=\frac{1}{2}+\frac{\epsilon}{2}$,
$b^{\prime}=-\frac{1}{2}+\epsilon,$ and $\sigma=\frac{1}{2}+\epsilon$.
Then for any
$u_{j}\in X_{0}^{s_{1},s_{2},b}(j=1,2)$,
we have
\begin{align}
      \|Q_{20}(u_{1},u_{2})\|_{X_{\sigma}^{s_{1}-3\epsilon,s_{2},
         b^{\prime}+\epsilon}}
\leq C\|u_{1}\|_{X_{0}^{s_{1},s_{2},b}}\,
      \|u_{2}\|_{X_{0}^{s_{1},s_{2},b}}.
      \label{3.035}
\end{align}
\end{Lemma}

Lemma 3.8 can be proved similarly to Lemma 4.10 of \cite{Hadac2008}.
\begin{Lemma}\label{Lemma3.9}
Let $s_{1}\geq-\frac{1}{2}+8\epsilon,s_{2}\geq0$,
$b=\frac{1}{2}+\frac{\epsilon}{2}$,
$b^{\prime}=-\frac{1}{2}+\epsilon,$ and
$\sigma=\frac{1}{2}+\epsilon$.
Then for any
$u_{j}\in X_{0}^{s_{1},s_{2},b}(j=1,2)$.
we have
\begin{align}
       \|Q_{20}(u_{1},u_{2})\|_{X_{\sigma}^{s_{1},
           s_{2}, b^{\prime}}}
\leq C \|u_{1}\|_{X_{0}^{s_{1},s_{2},b}}\,
       \|u_{2}\|_{X_{0}^{s_{1},s_{2},b}}.
       \label{3.036}
\end{align}
\end{Lemma}
{\bf Proof.}
By duality, it suffices to  prove
\begin{align}
      \left|\int_{\SR^{3}}\bar{u}Q_{20}(u_{1},u_{2})dxdydt\right|
\leq C\|u\|_{X_{-\sigma}}^{-s_{1} ,-s_{2},-b^{\prime} }\,
      \|u_{1}\|_{X_{0}^{s_{1},s_{2},b}}\,
      \|u_{2}\|_{X_{0}^{s_{1},s_{2},b}}.
\label{3.48}
\end{align}
for any $u\in X_{-\sigma}^{-s_{1},-s_{2},-b^{\prime}}.$

Let $F_{j}(j=1,2)$  be the same as in \eqref{Fj}, $G_2(\zeta)$
and $K_2(\zeta_1, \zeta)$ be the same as in \eqref{G2}
and \eqref{K2} respectively.
Then (\ref{3.48}) is equivalent to
\begin{align}
       \left|\int_{\SR^{6}}
      \chi_{_{A_{20}}}\,
      K_2(\zeta_1, \zeta)\,
      G_2(\zeta)\,
      F_{1}(\zeta_1) \,
      F_{2}(\zeta_2)\,d\zeta_1 d\zeta\right|
\leq C\|G_{2}\|_{L_{\zeta}^{2}}\,
      \|F_{1} \|_{L_{\zeta}^{2}}\,
      \|F_{2} \|_{L_{\zeta}^{2}}.
      \label{3.038}
\end{align}
We argue in two cases:  $|\xi|<1$ and $|\xi|\geq1$.

In the first case, $|\xi|<1$, note that $s_{1}\geq-\frac{1}{2}+8\epsilon,$
 we have
\begin{align*}
      K_2(\zeta_1, \zeta)
\leq  \frac{C|\xi_1|^{-2s_{1}-1+2\epsilon}}
      {\langle \lambda_1\rangle^{b}\,
      \langle \lambda_{2}\rangle^{b}}
\leq  \frac{C}
      {\langle \lambda_1\rangle^{b}\,
      \langle \lambda_{2}\rangle^{b}}.
\end{align*}
In the second case, $|\xi|\geq1$, note that
$s_{1}\geq-\frac{1}{2}+8\epsilon$ and $\sigma=\frac{1}{2}+\epsilon,$
 we have
\begin{align*}
      K_2(\zeta_1, \zeta)
\leq  \frac{C|\xi_1|^{-s_{1}+\sigma-1+2\epsilon}}
      {\langle \lambda_1\rangle^{b}\,
      \langle \lambda_{2}\rangle^{b}}
\leq  \frac{C|\xi_1|^{-5\epsilon}}
      {\langle \lambda_1\rangle^{b}\,
      \langle \lambda_{2}\rangle^{b}}
\leq  \frac{C}
      {\langle \lambda_1\rangle^{b}\,
      \langle \lambda_{2}\rangle^{b}}.
\end{align*}
Thus, for both cases, (\ref{3.038}) can be proved
 similarly to Part (i) of Lemma \ref{Lemma3.5}.

\begin{Lemma}\label{Lemma3.10}
Let $s_{1}\geq-\frac{1}{2}+8\epsilon,s_{2}\geq0$,
$b=\frac{1}{2}+\frac{\epsilon}{2}$,
$b^{\prime}=-\frac{1}{2}+\epsilon,$
and $\sigma=\frac{1}{2}+\epsilon$.
Then for any
$u_{j}\in X_{0}^{s_{1},s_{2},b}\,(j=1,2)$,
we have
\begin{align}
     \|Q_{21}(u_{1},u_{2})\|_{X_{\sigma}^{s_{1}-3\epsilon,
            s_{2}, b^{\prime}+\epsilon}}
\leq C\|u_{1}\|_{X_{0}^{s_{1},s_{2},b}}\,
      \|u_{2}\|_{X_{0}^{s_{1},s_{2},b}}.
      \label{3.039}
\end{align}
\end{Lemma}
{\bf Proof.}
By duality, it suffices to  prove
\begin{align}
      \left|\int_{\SR^{3}}\bar{u}\,
      Q_{21}(u_{1},u_{2})dxdydt\right|
\leq C\|u\|_{X_{-\sigma}^{-s_{1}+3\epsilon,
                -s_{2},-b^{\prime}-\epsilon}}\,
      \|u_{1}\|_{X_{0}^{s_{1},s_{2},b}}\,
      \|u_{2}\|_{X_{0}^{s_{1},s_{2},b}}
\label{3.51}
\end{align}
for any
$
 u\in X_{-\sigma}^{-s_{1}+3\epsilon,
      -s_{2},-b^{\prime}-\epsilon}.
$

Let $F_{j}(j=1,2)$  be the same as in \eqref{Fj},
$G_2(\zeta)$ be the same as in \eqref{G2}, and write
\begin{align}
    K_{5}(\zeta_1,\zeta)
=  \frac{|\xi|^{1-\sigma}\,
     \langle\xi\rangle^{s_{1}-3\epsilon+\sigma}\,
     \langle \eta\rangle ^{s_{2}}}
    {\langle \lambda\rangle^{-b^{\prime}-\epsilon}\,
     \langle \xi_{1}\rangle^{s_{1}}\,
     \langle \eta_{1}\rangle^{s_{2}}\,
     \langle \lambda_{1}\rangle^{b} \,
     \langle \xi_{2}\rangle^{s_{1}}\,
     \langle \eta_{2}\rangle^{s_{2}}\,
     \langle \lambda_{2}\rangle^{b}} .
\end{align}
Then (\ref{3.51}) is equivalent to
\begin{align}
       \left|\int_{\SR^{6}}
      \chi_{_{A_{21}}}
      K_5(\zeta_1, \zeta)\,
      G_2(\zeta)\,
      F_{1}(\zeta_1) \,
      F_{2}(\zeta_2)\,d\zeta_1 d\zeta\right|
\leq C\|G_{2}\|_{L_{\zeta}^{2}}\,
      \|F_{1} \|_{L_{\zeta}^{2}}\,
      \|F_{2} \|_{L_{\zeta}^{2}}.
      \label{3.041}
\end{align}
We argue in two cases:  $|\xi|<1$ and $|\xi|\geq1$.

In the first case, $|\xi|<1$, we consider
\begin{align}
&\left|1-\frac{1}{3\xi^{2}\xi_{2}^{2}}\right|\geq \frac{1}{4}, \label{3.054}\\
&\left|1-\frac{1}{3\xi^{2}\xi_{2}^{2}}\right|< \frac{1}{4}. \label{3.055}
\end{align}
Note that
$
s_{1}\geq-\frac{1}{2}+8\epsilon
$,
when (\ref{3.054}) is valid,
$
   \langle\lambda\rangle^{-b^{\prime}}
   \langle\lambda_{1}\rangle^{-b}
\leq
   \langle\lambda_{1}\rangle^{-b^{\prime}}
   \langle\lambda\rangle^{-b}
$,
we have
\begin{align*}
      K_{5}(\zeta_1,\zeta)
\leq \frac{C|\xi|^{1-\sigma}
     \langle\xi\rangle^{s_{1}+\sigma}}
    {\langle \lambda\rangle^{b}
     \langle \xi_{1}\rangle^{s_{1}}
     \langle \lambda_{1}\rangle^{-b^{\prime}}
     \langle \xi_{2}\rangle^{s_{1}}
     \langle \lambda_{2}\rangle^{b}}
\le  \frac{C|\xi_{1}|^{-2s_{1}-1+4\epsilon}}
    {\langle \lambda\rangle^{b}
     \langle \lambda_{2}\rangle^{b}}
\leq \frac{C}
     {\langle \lambda\rangle^{b}
      \langle \lambda_{2}\rangle^{b}}.
\end{align*}
When (\ref{3.055}) is valid,
$
   \langle\lambda\rangle^{-b^{\prime}-\epsilon}
   \langle\lambda_{1}\rangle^{-b}
\leq
   \langle\lambda_{1}\rangle^{-b^{\prime}-\epsilon}
   \langle\lambda\rangle^{-b}
$,
we have
\begin{align*}
      K_{5}(\zeta_1,\zeta)
\leq \frac{C|\xi|^{1-\sigma}
     \langle\xi\rangle^{s_{1}+\sigma-3\epsilon}}
    {\langle \lambda\rangle^{b}\,
     \langle \xi_{1}\rangle^{s_{1}}\,
     \langle \lambda_{1}\rangle^{-b'-\epsilon}\,
     \langle \xi_{2}\rangle^{s_{1}}\,
     \langle \lambda_{2}\rangle^{b}}
\le  \frac{C|\xi|^{-\epsilon}\,
     |\xi_{1}|^{-2s_{1}-1+4\epsilon}}
    {\langle \lambda\rangle^{b}\,
     \langle \lambda_{2}\rangle^{b}}
\leq \frac{C}
     {\langle \lambda\rangle^{b}\,
      \langle \lambda_{2}\rangle^{b}}.
\end{align*}
In the case  $|\xi|\geq1$, note that $s_{1}\geq-\frac{1}{2}+8\epsilon$
 and $\sigma=\frac{1}{2}+\epsilon,$
 we have
\begin{align*}
      K_{5}(\zeta_1,\zeta)
\leq \frac{C|\xi|^{\epsilon}
      |\xi_{1}|^{-s_{1}+\sigma-1+4\epsilon}}
     {\langle \lambda\rangle^{b}\,
      \langle \lambda_{2}\rangle^{b}}
\leq \frac{C|\xi_{1}|^{-s_{1}-1+\sigma+5\epsilon}}
     {\langle\lambda\rangle^{b}\,
      \langle\lambda_{2}\rangle^{b}}
\leq \frac{ C}
     {\langle \lambda\rangle^{b}\,
      \langle \lambda_{2}\rangle^{b}}.
\end{align*}
For all these cases, \eqref{3.041} can be proved similarly to Lemma 3.1.

This completes the proof of this lemma.

\begin{Lemma}\label{Lemma3.11}
Let
$s_{1}\geq-\frac{1}{2}+8\epsilon,s_{2}\geq0$,
$b=\frac{1}{2}+\frac{\epsilon}{2}$,
$b^{\prime}=-\frac{1}{2}+\epsilon,$ and $\sigma=\frac{1}{2}+\epsilon$.
Then for any $u_{j}\in X_{0}^{s_{1},s_{2},b}\,(j=1,2)$,
we have
\begin{align}
      \|Q_{21}(u_{1},u_{2})\|_{X_{\sigma}^{s_{1},s_{2},
             b^{\prime}}}
& \leq C\|u_{1}\|_{X_{0}^{s_{1},s_{2},b}}\,
      \|u_{2}\|_{X_{0}^{s_{1},s_{2},b}},
      \\ \|Q_{22}(u_{1},u_{2})\|_{X_{\sigma}^{s_{1}-3\epsilon,s_{2},
             b^{\prime}+\epsilon}}
&\leq C\|u_{1}\|_{X_{0}^{s_{1},s_{2},b}}\,
      \|u_{2}\|_{X_{0}^{s_{1},s_{2},b}}
      \\
      \|Q_{22}(u_{1},u_{2})\|_{X_{\sigma}^{s_{1},
       s_{2},b^{\prime}}}
& \leq C\|u_{1}\|_{X_{0}^{s_{1},s_{2},b}}\,
      \|u_{2}\|_{X_{0}^{s_{1},s_{2},b}}.
      \label{3.044}
\end{align}
\end{Lemma}
This lemma can be proved similarly to Lemma 3.10.

 Now we prove Theorem 1.1.

It is obvious that,
$
\partial_{x}(u_{1}u_{2})=\partial_{x}P_{1}(u_{1},u_{2})
+\partial_{x}P_{1}(u_{2},u_{1}).
$
Consequently, by symmetry, it suffices to prove
$
    \|\partial_{x}P_{1}(u_{1}u_{2})\|_{\widetilde{X}}
\leq C \|u_{1}\|_{X}\>\|u_{2}\|_{X}.
$
Thus, we have to  prove
\begin{align}
     \left\|Q_{kj}(u_{1},u_{2})\right\|_{\tilde{X}}
\leq C \|u_{1}\|_{X}\>\|u_{2}\|_{X},
     \quad 0\le k,j\le 2. \label{Qkj}
\end{align}
By  the definition of $\tilde{X}$ and Lemma 3.1
 as well as \eqref{embedding}, we have
\begin{align*}
     \|Q_{00}(u_{1},u_{2})\|_{\tilde{X}}
&\leq C\|Q_{00}(u_{1},u_{2})\|_{\tilde{X}_{1}}\\
&\leq C\|u_{1}\|_{X_{0}^{s_{1},s_{2},b}}\,
      \|u_{2}\|_{X_{0}^{s_{1},s_{2},b}}
\leq C\|u_{1}\|_{X}\,
      \|u_{2}\|_{X}.
\end{align*}
Combining Lemmas 3.2, 3.3, 3.7--3.11 with (\ref{embedding}),
we have
\begin{align}\nonumber
      \|Q_{kj}(u_{1},u_{2})\|_{\tilde{X}}
&\leq C\|Q_{kj}(u_{1},u_{2})\|_{\tilde{X}_{2}}\\
&\leq C\|u_{1}\|_{X_{0}^{s_{1},s_{2},b}}\,
      \|u_{2}\|_{X_{0}^{s_{1},s_{2},b}}
 \leq C\|u_{1}\|_{X}\,
      \|u_{2}\|_{X}\label{3.045}
\end{align}
for all rest $Q_{kj}$ except $Q_{11}.$ By using a proof
similar to  \cite{Hadac2008}, we can prove that \eqref{Qkj}
holds for $Q_{11}$ with the aid of Lemmas 3.4-3.6.

This completes the proof of Theorem  1.1.

\begin{Lemma}\label{Lemma3.12}
Let $A_{N}:=\left\{(\tau,\xi,\eta)|\frac{N}{2}\leq |\xi|\leq 2N\right\}$ and $u_{N_{1}},u_{N_{2}}$ possess the same support defined as in Lemma 2.18. Then,
there exists $C>0$ such that for all $T>0$ and functions $u_{N_{1}},u_{N_{2}},u_{N_{3}}\in V_{-,S}^{2}$ satisfies
$\mathscr{F}u_{N_{1}}\subset A_{N_{1}},$ $\mathscr{F}u_{N_{2}}\subset A_{N_{2}},$ and $\mathscr{F}w_{N_{3}}\subset A_{N_{3}}$ for dyadic numbers $N_{1},N_{2},N_{3}$,   the following inequality  holds true:
If $N_{2}\sim N_{3},$ then
\begin{eqnarray}
&&\left|\sum\limits_{N_{1}\leq  12 N_{2}}\int_{0}^{T}\int_{\SR^{2}}u_{N_{3}}\prod\limits_{j=1}^{2}u_{N_{j}}dxdydt\right|\nonumber\\&&
\leq C\left(\sum\limits_{N_{1}\leq 12 N_{2}}N_{1}^{-1}\|u_{N_{1}}\|_{V_{S}^{2}}^{2}\right)^{\frac{1}{2}}
\prod\limits_{j=2}^{3}N_{j}^{-\frac{1}{2}}\|u_{N_{j}}\|_{V_{S}^{2}},\label{3.061}
\end{eqnarray}
and if $N_{1}\sim N_{2},$ then
\begin{eqnarray}
&&\left(\sum\limits_{N_{3}\leq  \frac{1}{12 } N_{2}}N_{3}\sup\limits_{\|u_{N_{3}}\|_{V_{S}^{2}=1}}\left|\int_{0}^{T}\int_{\SR^{2}}u_{N_{3}}
\prod\limits_{j=1}^{2}u_{N_{j}}dxdydt\right|^{2}\right)^{\frac{1}{2}}\nonumber\\&&
\leq C\left(\sum\limits_{N_{3}\leq \frac{1}{12 }N_{2}}N_{1}^{-1}\|u_{N_{1}}\|_{V_{S}^{2}}^{2}\right)^{\frac{1}{2}}
\prod\limits_{j=1}^{2}N_{j}^{-\frac{1}{2}}\|u_{N_{j}}\|_{V_{S}^{2}},\label{3.062}
\end{eqnarray}
\end{Lemma}

Lemma 3.12 can be proved similarly to Proposition 3.1 of \cite{Hadac2008}.
\begin{Lemma}\label{Lemma3.13}
Let $A_{N}:=\left\{(\tau,\xi,\eta)|\frac{N}{2}\leq |\xi|\leq 2N\right\}$ and $u_{N_{1}},u_{N_{2}}$ possess the same support defined as in Lemma 2.19. Then,
there exists $C>0$ such that for all $T>0$ and functions $u_{N_{1}},u_{N_{2}},u_{N_{3}}\in V_{-,S}^{2}$ satisfies
$\mathscr{F}u_{N_{1}}\subset A_{N_{1}},$ $\mathscr{F}u_{N_{2}}\subset A_{N_{2}},$ and $\mathscr{F}w_{N_{3}}\subset A_{N_{3}}$ for dyadic numbers $N_{1},N_{2},N_{3}$,   the following holds true:
If $N_{2}\sim N_{3},$ then
\begin{eqnarray}
&&\left|\sum\limits_{N_{1}\leq  12 N_{2}}\int_{0}^{T}\int_{\SR^{2}}u_{N_{3}}\prod\limits_{j=1}^{2}u_{N_{j}}dxdydt\right|\nonumber\\&&
\leq C\left(\sum\limits_{N_{1}\leq 12 N_{2}}N_{1}^{-1}\|u_{N_{1}}\|_{V_{S}^{2}}^{2}\right)^{\frac{1}{2}}
\prod\limits_{j=2}^{3}N_{j}^{-\frac{1}{2}}\|u_{N_{j}}\|_{V_{S}^{2}},\label{3.063}
\end{eqnarray}
and if $N_{1}\sim N_{2},$ then
\begin{eqnarray}
&&\left(\sum\limits_{N_{3}\leq  \frac{1}{12 } N_{2}}N_{3}\sup\limits_{\|u_{N_{3}}\|_{V_{S}^{2}=1}}\left|\int_{0}^{T}\int_{\SR^{2}}u_{N_{3}}
\prod\limits_{j=1}^{2}u_{N_{j}}dxdydt\right|^{2}\right)^{\frac{1}{2}}\nonumber\\&&
\leq C\left(\sum\limits_{N_{3}\leq \frac{1}{12 }N_{2}}N_{1}^{-1}\|u_{N_{1}}\|_{V_{S}^{2}}^{2}\right)^{\frac{1}{2}}
\prod\limits_{j=1}^{2}N_{j}^{-\frac{1}{2}}\|u_{N_{j}}\|_{V_{S}^{2}}.\label{3.064}
\end{eqnarray}
\end{Lemma}

Lemma 3.13 can be proved similarly to Proposition 3.1 of \cite{Hadac2008}.
\begin{Lemma}\label{Lemma3.14}
Let $A_{N}:=\left\{(\tau,\xi,\eta)|\frac{N}{2}\leq |\xi|\leq 2N\right\}$ and $u_{N_{1}},u_{N_{2}}$ possess the same support defined as in Lemma 2.20. Then,
there exists $C>0$ such that for all $T>0$ and functions $u_{N_{1}},u_{N_{2}},u_{N_{3}}\in V_{-,S}^{2}$ satisfies
$\mathscr{F}u_{N_{1}}\subset A_{N_{1}},$ $\mathscr{F}u_{N_{2}}\subset A_{N_{2}},$ and $\mathscr{F}w_{N_{3}}\subset A_{N_{3}}$ for dyadic numbers $N_{1},N_{2},N_{3}$,   the following inequality holds true:
If $N_{2}\sim N_{3},$ then
\begin{eqnarray}
&&\left|\sum\limits_{N_{1}\leq  12 N_{2}^{-1}}\int_{0}^{T}\int_{\SR^{2}}u_{N_{3}}\prod\limits_{j=1}^{2}u_{N_{j}}dxdydt\right|\nonumber\\&&
\leq C\left(\sum\limits_{N_{1}\leq 12 N_{2}^{-1}}N_{1}^{-1}\|u_{N_{1}}\|_{V_{S}^{2}}^{2}\right)^{\frac{1}{2}}
\prod\limits_{j=2}^{3}N_{j}^{-\frac{1}{2}}\|u_{N_{j}}\|_{V_{S}^{2}},\label{3.065}
\end{eqnarray}
and if $N_{1}\sim N_{2},$ then
\begin{eqnarray}
\left(\sum\limits_{N_{3}\leq  12 N_{2}^{-1}}N_{3}\sup\limits_{\|u_{N_{3}}\|_{V_{S}^{2}=1}}\left|\int_{0}^{T}\int_{\SR^{2}}u_{N_{3}}
\prod\limits_{j=1}^{2}u_{N_{j}}dxdydt\right|^{2}\right)^{\frac{1}{2}}
\leq C
\prod\limits_{j=1}^{2}N_{j}^{-\frac{1}{2}}\|u_{N_{j}}\|_{V_{S}^{2}}.\label{3.066}
\end{eqnarray}
\end{Lemma}
\noindent{\bf Proof.} We define $\tilde{u}_{N_{j}}=\chi_{[0,T)}u_{N_{j}}(j=1,2,3)$ and $M=\frac{N_{1}N_2N_{3}}{8}$.
We decompose $Id=Q_{<M}+Q_{\geq M}$. Then, the left-hand side of (\ref{3.065}) can be decomposed into eight pieces
\begin{eqnarray}
\int_{\SR}\prod\limits_{j=1}^{3}Q_{j}^{S}\tilde{u}_{N_{j}}dxdydt\label{3.067}
\end{eqnarray}
with $Q_{i}^{S}\in \left\{Q_{\geq M}^{S},Q_{<M}^{S}\right\}(i=1,2,3).$
From the analysis of Proposition 3.1 of \cite{Hadac2008}, we only consider the case that
$Q_{i}^{S}=Q_{\geq M}^{S}$ for some $1\leq i\leq 3.$

\noindent Case $Q_{1}^{S}=Q_{\geq M}^{S}$ can be proved similarly to case  $Q_{1}^{S}=Q_{\geq M}^{S}$ of Proposition
3.1 of \cite{Hadac2008}.

Now we consider case $Q_{2}^{S}=Q_{\geq M}^{S}$. Combining (\ref{2.031}) with (\ref{2.051}), by using the Cauchy-Schwarz  inequality with respect to $N_{1},$ we have that
\begin{eqnarray}
&&\left|\int_{\SR^{3}}Q_{1}^{S}\tilde{u}_{N_{1}}Q_{\geq M}^{S}\tilde{u}_{N_{1}}Q_{3}^{S}\tilde{u}_{N_{3}}dxdydt\right|
\nonumber\\&&\leq \left\|Q_{\geq M}^{S}\tilde{u}_{N_{2}}\right\|_{L^{2}(\SR^{3})}\left(\frac{1}{N_{3}}\right)^{\frac{15}{28}}
\left\|Q_{1}^{S}\tilde{u}_{N_{1}}\right\|_{V_{S}^{2}}\left\|Q_{3}^{S}\tilde{u}_{N_{3}}\right\|_{V_{S}^{2}}\nonumber\\
&&\leq \frac{C}{(N_{1}N_{2}N_{3})^{\frac{1}{2}}}\|u_{N_{2}}\|_{V_{S}^{2}}\left(\frac{1}{N_{3}}\right)^{\frac{15}{28}}
\left\|\tilde{u}_{N_{1}}\right\|_{V_{S}^{2}}\left\|\tilde{u}_{N_{3}}\right\|_{V_{S}^{2}}\nonumber\\
&&\leq C\left(\sum\limits_{N_{1} \leq 12N_{2}^{-1}}N_{1}^{-1}\|u_{N_{1}}\|_{V_{S}^{2}}^{2}\right)^{\frac{1}{2}}
\left(\sum\limits_{N_{1} \leq 12N_{2}^{-1}}\right)^{\frac{1}{2}}\left(\frac{1}{N_{3}}\right)^{\frac{15}{28}}
\prod\limits_{j=2}^{3}N_{j}^{-\frac{1}{2}}\|u_{N_{j}}\|_{V_{S}^{2}}\nonumber\\
&&\leq C\left(\sum\limits_{N_{1} \leq 12N_{2}^{-1}}N_{1}^{-1}\|u_{N_{1}}\|_{V_{S}^{2}}^{2}\right)^{\frac{1}{2}}
\prod\limits_{j=2}^{3}N_{j}^{-\frac{1}{2}}\|u_{N_{j}}\|_{V_{S}^{2}}.\label{3.068}
\end{eqnarray}
(\ref{3.066}) can be proved similarly to (38) of Proposition 3.1 of \cite{Hadac2008}.

This completes the proof of Lemma 3.14.

\begin{Lemma}\label{Lemma3.15}
Let $A_{N}:=\left\{(\tau,\xi,\eta)|\frac{N}{2}\leq |\xi|\leq 2N\right\}$ and $u_{N_{1}},u_{N_{2}}$ possess the same support defined as in Lemma 2.21. Then,
there exists $C>0$ such that for all $T>0$ and functions $u_{N_{1}},u_{N_{2}},u_{N_{3}}\in V_{-,S}^{2}$ satisfies
$\mathscr{F}u_{N_{1}}\subset A_{N_{1}},$ $\mathscr{F}u_{N_{2}}\subset A_{N_{2}},$ and $\mathscr{F}w_{N_{3}}\subset A_{N_{3}}$ for dyadic numbers $N_{1},N_{2},N_{3}$,   the following inequality holds true:
If $N_{2}\sim N_{3},$ then
\begin{eqnarray}
&&\left|\sum\limits_{N_{1}\leq  12 N_{2}^{-1}}\int_{0}^{T}\int_{\SR^{2}}u_{N_{3}}\prod\limits_{j=1}^{2}u_{N_{j}}dxdydt\right|\nonumber\\&&
\leq C\left(\sum\limits_{N_{1}\leq 12 N_{2}^{-1}}N_{1}^{-1}\|u_{N_{1}}\|_{V_{S}^{2}}^{2}\right)^{\frac{1}{2}}
\prod\limits_{j=2}^{3}N_{j}^{-\frac{1}{2}}\|u_{N_{j}}\|_{V_{S}^{2}},\label{3.069}
\end{eqnarray}
and if $N_{1}\sim N_{2},$ then
\begin{eqnarray}
\left(\sum\limits_{N_{3}\leq  12 N_{2}^{-1}}N_{3}\sup\limits_{\|u_{N_{3}}\|_{V_{S}^{2}=1}}\left|\int_{0}^{T}\int_{\SR^{2}}u_{N_{3}}
\prod\limits_{j=1}^{2}u_{N_{j}}dxdydt\right|^{2}\right)^{\frac{1}{2}}
\leq C
\prod\limits_{j=1}^{2}N_{j}^{-\frac{1}{2}}\|u_{N_{j}}\|_{V_{S}^{2}}.\label{3.070}
\end{eqnarray}
\end{Lemma}

Combining Lemma 2.21 with a proof similar to Lemma 3.14, we have that
Lemma 3.15 is valid.

\begin{Lemma}\label{Lemma3.16}
For all $0<T<\infty$ and for all $u_{1},u_{2}\in \dot{Z}^{-\frac{1}{2}}
\cap  C(\R;H^{1,1}(\R^{2}))$, we have that
\begin{eqnarray}
&&\left\|I_{T}(u_{1},u_{2})\right\|_{\dot{Z}^{-\frac{1}{2}}}\leq
C\prod\limits_{j=1}^{2}\|u_{j}\|_{ \dot{Y}^{-\frac{1}{2}}},\label{3.071}\\
&&\left\|I_{T}(u_{1},u_{2})\right\|_{\dot{Z}^{-\frac{1}{2}}}\leq C
\prod\limits_{j=1}^{2}\|u_{j}\|_{ \dot{Z}^{-\frac{1}{2}}},\label{3.072}
\end{eqnarray}
where
\begin{eqnarray*}
I_{T}(u_{1},u_{2})=\int_{0}^{t}\chi_{[0,T)}W(t-t^{\prime})
\partial_{x}(u_{1}u_{2})(t^{\prime})dt^{\prime}.
\end{eqnarray*}

\end{Lemma}
\noindent {\bf Proof.} Combining Lemmas 3.12-3.15 with Theorem 3.2 of
 \cite{Hadac2008}, we have that
(\ref{3.071}) is valid. Since $\dot{Z}^{-\frac{1}{2}}\subset
 \dot{Y}^{-\frac{1}{2}},$ from (\ref{3.071}),
we have that (\ref{3.072}) is valid.

This completes the proof of Lemma 3.16.

\begin{Lemma}\label{Lemma3.17}
For $0<T\leq1$ and $\epsilon>0$  and $u_{N_{1}}\in X,u_{N_{2}}\in
U_{S}^{2},u_{N_{3}}\in V_{-,S}^{2}$ with $\mathscr{F}u_{N_{j}}\subset
 A_{N_{j}}(j=1,2,3)$ for dyadic numbers $N_{j}(j=1,2,3)$,
we have
\begin{eqnarray}
\left|\int_{0}^{T}\int_{\SR^{2}}\prod\limits_{j=1}^{3}u_{N_{j}}dxdydt\right|\leq \frac{C(TN_{1})^{\frac{1}{4}-\epsilon}}{(N_{2}N_{3})^{\frac{1}{2}}}\|u_{N_{1}}\|_{X}\|u_{N_{2}}\|_{U_{S}^{2}}
\|u_{N_{3}}\|_{V_{S}^{2}},\label{3.073}
\end{eqnarray}
where $N_{1}\leq1 \leq N_{2}.$
\end{Lemma}
\noindent {\bf Proof.} Combining Lemmas 2.18-2.21 with Proposition 3.5 of \cite{Hadac2008}, we have that
Lemma 3.17 is valid.

This completes the proof of Lemma 3.17.

\begin{Lemma}\label{Lemma3.18}
For all $0<T<\infty$ and for all $u_{1},u_{2}\in \dot{Z}^{-\frac{1}{2}}
\cap  C(\R;H^{1,1}(\R^{2}))$, we have that
\begin{eqnarray}
&&\left\|I_{1}(u_{1},u_{2})\right\|_{\dot{Z}^{-\frac{1}{2}}}\leq
C\prod\limits_{j=1}^{2}\|u_{j}\|_{ Z^{-\frac{1}{2}}},\label{3.074}\\
&&\left\|I_{1}(u_{1},u_{2})\right\|_{Z^{-\frac{1}{2}}}\leq C
\prod\limits_{j=1}^{2}\|u_{j}\|_{ Z^{-\frac{1}{2}}},\label{3.075}
\end{eqnarray}
where
\begin{eqnarray*}
I_{T}(u_{1},u_{2})=\int_{0}^{t}\chi_{[0,T)}W(t-t^{\prime})
\partial_{x}(u_{1}u_{2})(t^{\prime})dt^{\prime}.
\end{eqnarray*}

\end{Lemma}
\noindent {\bf Proof.} Combining Propositions 1.6, 1.8, Lemmas 3.16-3.17 with a proof
similar to Theorem 3.6 of \cite{Hadac2008}, we have that(\ref{3.073})  is valid.
By using $\dot{Z}^{-\frac{1}{2}}\subset Z^{-\frac{1}{2}},$ we have that (\ref{3.074})
is valid.

This completes the proof of Lemma 3.18.

\bigskip

\bigskip

\noindent {\large\bf 4. Proof of Theorem  1.2}

\setcounter{equation}{0}

 \setcounter{Theorem}{0}

\setcounter{Lemma}{0}

\setcounter{section}{4}

In this section, combining Lemmas 2.3, Theorem 1.1 with the
fixed point theorem, we present the proof of Theorem 1.2.

Let  $0<T\leq 1,$ and $u_{1},u_{2}$ be
rapidly decreasing functions.
We define the bilinear operator $\Gamma_{T}$ by (\ref{4.03}).
We define
\begin{align}
   \Gamma_{T}(u_{1},u_{2})
=  \psi\left(\frac{t}{T}\right)
   \int_{0}^{t}W(t-\tau)\partial_{x}(u_{1}u_{2})d\tau.
   \label{4.01}
\end{align}
When
 $s_{1}>-\frac{1}{2},\,s_{2}\geq0$,
we can choose $b=\frac{1}{2}+\frac{\epsilon}{2}$,
 $b^{\prime}=-\frac{1}{2}+\epsilon$,
$\sigma=\frac{1}{2}+\epsilon$
so that Lemma 2.3 applies and we get
\begin{align}
&      \left\|\psi W(t)u_{0}\right\|_{X}
 \leq C\|u\|_{H^{s_{1},s_{2}}(\SR^{2})},
       \label{4.02}\\
&      \left\|\psi\Big(\frac{t}{T}\Big)
       \int_{0}^{t}W(t-\tau)Fd\tau\right\|_{X}
\leq  CT^{1-b+b^{\prime}}\|F\|_{\tilde{X}}.
       \label{4.03}
\end{align}
By the bilinear estimate in Theorem 1.1,
(\ref{4.01}) and (\ref{4.03}), we have
\begin{align}
     \left\|\Gamma_{T}(u_{1},u_{2})\right\|_{X}
\leq CT^{1+b^{\prime}-b}
     \left\|\partial_{x}(u_{1},u_{2})\right\|_{\tilde{X}}
\leq CT^{1+b^{\prime}-b}
     \|u_{1}\|_{X}\,
     \|u_{2}\|_{X},\label{4.04}
\end{align}
 Thus, we can extend $\Gamma_{T}$ to a continuous bilinear operator
  $\Gamma_{T}:X \times X\longrightarrow X.$
  As $\Gamma_{T}(u_{1},u_{2}) |_{[-T,T]}$
only depends on $u_{j}|_{[-T,T]}(j=1,2)$, $\Gamma_{T}$ also defines
 a continuous  bilinear operator $\Gamma_{T}$  from $X_{T}\times X_{T}$
 to $ X_{T}$.

 We define
\begin{align}
&\Phi_{T}(u,u_{0})=\psi W(t)u_{0}+\frac{1}{2}\Gamma_{T}(u,u),
u\in X_{T}, u_{0}
\in H^{s_{1},s_{2}}(\R^{2}),\label{4.05}\\
&B_{R}:=\left\{u_{0}\in H^{s_{1},s_{2}}(\R^{2})\mid
\|u_{0}\|_{H^{s_{1},s_{2}}}
\leq R\right\},\\
&A_{r}:=\left\{u\in X_{T}\mid \|u\|_{X_{T}}\leq r\right\}.
\label{4.06}
\end{align}
For $u_{0}\in B_{R},u\in A_{r},$ by using (\ref{4.02})--(\ref{4.05}),
we have
\begin{align}
      \left\|\Phi_{T}(u,u_{0})\right\|_{X}
&\leq \left\|\psi(t)W(t)u_{0}\right\|_{X}
     +\left\|\frac{1}{2}\psi\left(\frac{t}{\tau}\right)
      \int_{0}^{t}W(t-\tau)\partial_{x}(u^{2})d\tau\right\|_{X}\nonumber\\
&\leq C\|u_{0}\|_{H^{s_{1},s_{2}}}+CT^{1+b^{\prime}-b}
      \left\|\partial_{x}(u^{2})\right\|_{\tilde{X}}\nonumber\\
&\leq C\|u_{0}\|_{H^{s_{1},s_{2}}}+CT^{1+b^{\prime}-b}
      \left\|u\right\|_{X}^{2}\nonumber\\
&\leq CR+CT^{1+b^{\prime}-b}r^{2}.
      \label{4.08}
\end{align}
Here $C>0$ does not depend on $R,r$ and $T.$
Note that $1+b^{\prime}-b>0$,
we can choose $T\in (0,1)$ small such that
\begin{align}
T^{1+b^{\prime}-b}=\left[16C^{2}(R+1)\right]^{-1}.\label{4.09}
\end{align}
Combining (\ref{4.08}) with (\ref{4.09}), for $r=2CR,$ we have that
\begin{align}
&\left\|\Phi_{T}(u,u_{0})\right\|_{X}
\leq r.\label{4.010}
\end{align}
Thus, for fixed $u_{0}\in B_{R},$  $\Phi_{T}(\cdot,u_{0})$ maps
 $A_{r}$ into $A_{r}$.
By Theorem 1.1, (\ref{4.05})--(\ref{4.010}), we have
\begin{align*}
       \|\Phi_{T}(u, u_{0}) \|_{X}
&\leq C\left\|\Gamma_{T}(u-v,u+v)\right\|_{X}
 \leq CT^{1+b^{\prime}-b}\left\|u-v\right\|_{X}
       \left[\left\|u\right\|_{X}
        +\left\|v\right\|_{X}\right]
       \nonumber \\
&\leq  2CT^{1+b^{\prime}-b}r\left\|u-v\right\|_{X}
\leq   \frac{1}{2} \left\|u-v\right\|_{X}.
\end{align*}
Thus, $\Phi_{T}(\cdot,u_{0})$ is a contraction on $A_r$.
By the Banach fixed point theorem, $\Phi_{T}$ admits a
unique fixed point in $A_{r}$. The rest of the proof of
Theorem 1.2 is similar to \cite{Hadac2008}.

\bigskip

\bigskip

\noindent {\large\bf 5. Proof of Theorem 1.3}

\setcounter{equation}{0}

 \setcounter{Theorem}{0}

\setcounter{Lemma}{0}

\setcounter{section}{5}
Similar to \cite[Theorem 4.2]{KZ}, we prove Theorem 1.3.
The solution to the Cauchy problem for (\ref{1.04}) is
formally equivalent to the following integral equation
\begin{align}
u(u)=W(t)u_{0}
   -\frac{1}{2}\int_{0}^{t}W(t-\tau)\partial_{x}(u^{2})d\tau.
\label{5.01}
\end{align}
To prove Theorem 1.3, motivated by \cite[Lemma 4.2]{KZ},
it suffices to prove Lemma 5.1.
\begin{Lemma}\label{Lemma5.1}
Let $u(t)$ be the solution to  (\ref{5.01}). Then, if $s_{1}<-\frac{1}{2}$,
 then there is no $T>0$ for which the map
$u_{0}\rightarrow u(t),$ $t\in [0,T]$ is $C^{3}$ at zero from
 $H^{s_{1},0}(\R^{2})$ to $H^{s_{1},0}(\R^{2})$.
More precisely, for fixed point $t_{0}\geq0$ and $s_{1}<-\frac{1}{2}$,
 then the following inequality is invalid:
\begin{align}
\left\|\int_{0}^{t_{0}}W(t-\tau)\partial_{x}(W(\tau)u_{0} W_{2}(\tau))
\right\|_{H^{s_{1},0}}
\leq C\|u_{0}\|_{H^{s_{1},0}}^{3},\label{5.02}
\end{align}
where
$\displaystyle
W_{2}(x,y,\tau)=\frac{1}{2}\int_{0}^{\tau}W(\tau-t^{\prime})
\partial_{x}(W(t^{\prime})u_{0})\,dt^{\prime}.\label{5.03}
$
\end{Lemma}
{\bf Proof.} For  fixed $t_{0}<\frac{1}{200}$ and we define
 $\theta \in C_{0}^{\infty}(\R)$ by
$
\theta(t)=1
$
if $\frac{1}{2}t_{0}\leq t\leq 2t_{0}$
and
$
\theta(t)=0
$
if $0\leq t_{0}<\frac{1}{4}t_{0}$ or $ t\geq 4t_{0}.$
Thus, if $s_{1}<-\frac{1}{2}$, to show (\ref{5.02}),
 it suffices to prove that
\begin{align}
   \Big\|\theta(t)
   \int_{0}^{t_{0}}W(t-\tau)\partial_{x}
   (W(\tau)u_{0}W_{2}(\tau))d\tau\Big\|_{H^{s_{1},0}}
\leq C\|u_{0}\|_{H^{s_{1},0}}^{3} \label{5.04}
\end{align}
is invalid. Recall that
$\phi(\xi,\eta)=\xi^{3}-\frac{\eta^{2}+1}{\xi}$. We let
\begin{align*}
&  \chi_{_{+}}(\xi)
=  \chi_{_{[2^{k}+2^{1- {k}/{2}},\,2^{k}+2^{2-{k}/{2}}]}}(\xi),\quad
   \chi_{_{-}}(\xi)
=  \chi_{_{[-2^{k},\,-2^{k}+2^{-{k}/{2}}]}}(\xi),\\
&  \chi_{2} (\frac{\eta}{\xi})
=  \chi_{_{[2 ,2 +2^{-{k}/{2}}]}}\big({\eta}/{\xi}\big),\\
& P=P(\xi_{1},\xi_{2},\xi,\eta_{1},\eta_2,\eta)\\
&\phantom{P }= \phi(\xi-\xi_{1},\eta-\eta_{1})
  +\phi(\xi_{1}-\xi_{2},\eta_{1}-\eta_{2})
  +\phi(\xi_{2},\eta_{2})
  -\phi(\xi,\eta), \\
& Q(\xi_{1},\xi,\eta_{1},\eta)
= \phi(\xi_{1},\eta_{1})
 +\phi(\xi-\xi_{1},\eta-\eta_{1})
 -\phi(\xi,\eta).
\end{align*}
Choose $u_0$ via
$$
  \mathscr{F}_{xy}u_0
=\bigl(\chi_{_{+}}(\xi)
      +\chi_{_{-}}(\xi)\bigr)
 \chi_2(\frac{\eta}{\xi}).
$$
Similar to \cite[Lemma 4.2]{KZ}, we have
\begin{align}
&\mathscr{F}_{xy}
   \left[\theta(t)
    \int_{0}^{t_{0}}W(t-\tau)\partial_{x}
    (W(\tau)u_{0} W_{2}(\tau))d\tau\right]
    \nonumber\\
=&\theta(t)\xi e^{it_{0}\phi(\xi,\eta)}
   \int_{\SR^{4}}
   \left(
    \frac{e^{i t_{0}P}-1}
         {PQ(\xi_{2},\xi_{1},\eta_{2},\eta_{1})}
   -\frac{e^{i t_{0}Q(\xi_{1},\xi,\eta_{1},\eta)}-1}
         {Q(\xi_{1},\xi,\eta_{1},\eta)
          Q(\xi_{2},\xi_{1},\eta_{2},\eta_{1})}
    \right)\nonumber\\
&\qquad\times
    \xi_{1}\mathscr{F}_{xy}
     u_{0}(\xi-\xi_{1},\eta-\eta_{1})
     \mathscr{F}_{xy}u_0(\xi_{1}-\xi_{2},\eta_{1}-\eta_{2})
     d\xi_1 d\xi_2 d\eta_1 d\eta_2 ,
       \label{5.06}
\end{align}
Note that
\begin{align*}
 &  \phi(\xi_{1},\mu_{1})+\phi(\xi_{2},\mu_{2})
   -\phi(\xi_{1}+\xi_{2},\mu_{1}+\mu_{2})\nonumber\\
=&  -3\xi_{1}\xi_{2}(\xi_{1}+\xi_{2})
    - \frac{\xi_{1}\xi_{2}}{\xi_{1}+\xi_{2}}
      \left(\frac{\mu_{1}}{\xi_{1}}
           -\frac{\mu_{2}}{\xi_{2}}\right)^{2}
    - \frac{\xi_{1}^{2}-\xi_{1}\xi_{2}+\xi_{2}^{2}}
           {\xi_{1}\xi_{2}(\xi_{1}+\xi_{2})},
\end{align*}
we have
$$
     \left|Q(\xi_{1},\xi,\eta_{1},\eta)\right|
\sim 2^{{3k}/{2}},\quad
     \left|Q(\xi_{2},\xi_{1},\eta_{2},\eta_{1})\right|
\sim 2^{{3k}/{2}},
$$
and thus
\begin{align}
     \left|\frac{e^{i t_{0}P}-1}
     {PQ(\xi_{2},\xi_{1},\eta_{2},\eta_{1})}\right|
\sim 2^{- 3k/2 }t_{0},\quad
     \left|\frac{e^{i t_{0}Q(\xi_{1},\xi,\eta_{1},\eta)}-1}
     {Q(\xi_{1},\xi,\eta_{1},\eta)
     Q(\xi_{2},\xi_{1},\eta_{2},\eta_{1})}\right|
\sim C2^{- 3k }t_{0}.
\label{5.07}
\end{align}
Hence, similar to \cite[Lemma 4.2]{KZ}  we get
\begin{align*}
&   \mathscr{F}_{xy}
    \left[\theta(t)
     \int_{0}^{t_{0}}
     W(t-\tau)\partial_{x}(W(\tau)u_{0}(x,y)
     W_{2}(x,y,\tau))d\tau\right] \\
\sim &
     t_{0}2^{-k}\!\!
     \int_{\SR^{4}}\!\!
     \chi_{_{+}}(\xi-\xi_{1})
     \chi_{_{-}}(\xi_{1}-\xi_{2})\chi_{_{+}}(\xi_{2})
     \chi_{_{2}}\Big(\frac{\mu-\mu_{1}}{\xi-\xi_{1}}\Big)
     \chi_{_{2}}\Big(\frac{\mu_{1}-\mu_{2}}{\xi_{1}-\xi_{2}}\Big)
     \chi_{_{2}}\Big(\frac{\mu_{2}}{\xi_{2}}\Big)
     d\xi_{1}d\xi_{2}d\mu_{1}d\mu_{2},
\end{align*}
and therefore, we obtain
\begin{align*}
\left\|\theta(t)\int_{0}^{t_{0}}W(t-\tau)
\partial_{x}(W(\tau)u_{0}(x,y)W_{2}(x,y,\tau))
d\tau\right\|_{H^{s_{1},0}}\leq
C\|u_{0}\|_{H^{s_{1},0}}^{3}.
\end{align*}
From the above estimate we have
$$
   2^{-k}t_{0}\leq C2^{2ks_{1}}.
$$
When $k\to + \infty$, this implies that
$s_{1}\geq-\frac{1}{2}$,  we derive a contradiction.

The proof of Lemma 5.1 is completed.

\bigskip
\bigskip
\noindent {\large\bf 6. Proof of Theorem 1.4}

\setcounter{equation}{0}

 \setcounter{Theorem}{0}

\setcounter{Lemma}{0}

\setcounter{section}{6}
Inspired by Theorem 1.2 of \cite{Hadac2009}, we present the proof of Theorem 1.4.

We define
\begin{eqnarray}
B_{\lambda,R}=\left\{u_{0}\in H^{-\frac{1}{2},0}(\R^{2})|u_{0}=v_{0}+w_{0},\|v_{0}\|_{\dot{H}^{-\frac{1}{2},0}}<\lambda,\|w_{0}\|_{L^{2}}\leq R\right\},\label{6.01}
\end{eqnarray}
where $\lambda,R$ will be determined later.
Let $u_{0}\in B_{\lambda,R}$   with  $u_{0}=v_{0}+w_{0}$.
By using a direct computation, we have that
\begin{eqnarray}
&&\left\|\chi e^{tS}w_{0}\right\|_{X}\leq
 \left\|\chi e^{tS}w_{0}\right\|_{\dot{X}^{0,1,1}}\leq C\|w_{0}\|_{L^{2}}\leq CR,\label{6.02}\\
&&\left\|\chi e^{tS}v_{0}\right\|_{\dot{Z}^{-\frac{1}{2}}}\leq \left\|\chi e^{tS}v_{0}\right\|_{\dot{X}^{-\frac{1}{2},\frac{1}{2},1}}\leq C
\|v_{0}\|_{\dot{H}^{-\frac{1}{2}}}\leq C\lambda\label{6.03}.
\end{eqnarray}
Combining (\ref{6.02}) with (\ref{6.03}), we have that
\begin{eqnarray}
\left\|e^{tS}u_{0}\right\|_{{Z}^{-\frac{1}{2}}}\leq C(\lambda+R).\label{6.04}
\end{eqnarray}
We define
\begin{eqnarray}
D_{r}=\left\{u\in Z^{-\frac{1}{2}}([0,1])|\|u\|_{Z^{-\frac{1}{2}}([0,1])}
\leq r:=\frac{1}{8C+8}\right\}\label{6.05}.
\end{eqnarray}
Here, $C$ is from (\ref{3.075}).
We choose $R=\lambda=\frac{1}{(8C+8)^{2}}.$
We define
\begin{eqnarray}
\Phi(u)=e^{tS}u_{0}-\frac{1}{2}I_{1}(u,u)(t),t\in [0,1].\label{6.06}
\end{eqnarray}
By using (\ref{3.075}) and (\ref{6.02})-(\ref{6.06}), we have that
\begin{eqnarray}
\left\|\Phi(u)\right\|_{{Z}^{-\frac{1}{2}}}\leq \left\|e^{tS}u_{0}\right\|_{{Z}^{-\frac{1}{2}}}+\left\|\frac{1}{2}I_{1}(u,u)(t)\right\|\leq C
(\lambda+R)+Cr^{2}\leq \frac{r}{2}+\frac{r}{2}=r\label{6.07}
\end{eqnarray}
and
\begin{eqnarray}
&&\left\|\Phi(u)-\Phi(v)\right\|_{{Z}^{-\frac{1}{2}}}\leq \left\|\frac{1}{2}I_{1}(u,u)(t)\right\|_{{Z}^{-\frac{1}{2}}}\nonumber\\&&\leq C\left[\|u\|_{Z^{-\frac{1}{2}}([0,1])}+\|v\|_{Z^{-\frac{1}{2}}([0,1])}\right]\left\|u-v\right\|_{{Z}^{-\frac{1}{2}}}\leq 2Cr\left\|u-v\right\|_{{Z}^{-\frac{1}{2}}}\nonumber\\&&\leq \frac{1}{2}\left\|u-v\right\|_{{Z}^{-\frac{1}{2}}}.\label{6.08}
\end{eqnarray}
Thus, $\Phi$ is a strict contraction. Consequently,  $\Phi$ has a unique fixed point in $D_{\lambda}$, which solves (\ref{6.06}) on the interval $(0,1)$.
By using the implicit function theorem, we know that  the map  $B_{\lambda,R}\mapsto D_{r}$, $u_{0}\mapsto u$ is analytic. We According to the definition of $Z^{-\frac{1}{2}}([0,1])$, we also have the embedding $Z^{-\frac{1}{2}}([0,1])\hookrightarrow C([0,1];H^{-\frac{1}{2}}(\R^{2}))$.
 Now, we assume that $u_{0}\in B_{\lambda,R}$ with $R\geq \lambda=\frac{1}{(8C+8)^{2}}$. We define
 $u_{0\Lambda}:=\Lambda^{2}u_{0}\left(\Lambda x,\Lambda^{2}y\right)$ and $\Lambda=R^{-2}\lambda^{2}$.  Consequently£¬  we have that $u_{\Lambda}=\Lambda^{2}u\left(\Lambda x,\Lambda^{2}y, \Lambda^{3}t\right)$ which is the solution to the problem with the initial data $u_{\Lambda}(0)=u_{0\Lambda}$. In particular, $u_{\Lambda} \in Z^{-\frac{1}{2}}([0,1])$. By rescaling, we
 we have that $u\in Z^{-\frac{1}{2}}([0,\Lambda ^{6}R^{-6}])$ which is the solution to the problem with the initial data $u(0)=u_{0}$.

Now we claim that if $u_{1}(0)=u_{2}(0)$ and  $u_{1},u_{2}\in Z^{-\frac{1}{2}}([0,T])$ are solutions, then we have that
$u_{1}=u_{2}$ on $Z^{-\frac{1}{2}}([0,T])$. We decompose $u_{j}=v_{j}+w_{j},v_{j}\in X([0,T]),w_{j}\in \dot{Z}^{-\frac{1}{2}}([0,T])$, $w_{j}(0)=0$.
For small $0<\tau\leq T$, by using (\ref{3.075}), we have that
\begin{eqnarray}
&&\left\|u_{1}-u_{2}\right\|_{{Z}^{-\frac{1}{2}}([0,\tau])}\leq \left\|\frac{1}{2}I_{1}(u,u)(t)\right\|_{{Z}^{-\frac{1}{2}}([0,\tau])}\nonumber\\&&\leq C\tau^{\frac{1}{4}-\epsilon}\left(
\|v_{1}\|_{{Z}^{-\frac{1}{2}}([0,\tau])}+\|v_{2}\|_{{Z}^{-\frac{1}{2}}([0,\tau])}\right)
\left\|u_{1}-u_{2}\right\|_{{Z}^{-\frac{1}{2}}([0,\tau])}\nonumber\\&&\qquad+C\left(
\|w_{1}\|_{\dot{Z}^{-\frac{1}{2}}([0,\tau])}+\|v_{2}\|_{\dot{Z}^{-\frac{1}{2}}([0,\tau])}\right)
\left\|u_{1}-u_{2}\right\|_{{Z}^{-\frac{1}{2}}([0,\tau])([0,\tau])}\label{6.09}.
\end{eqnarray}
By using (ii) of Proposition 1.15, we have that
\begin{eqnarray}
\|w_{j}\|_{\dot{Z}^{-\frac{1}{2}}([0,\tau])}\leq \epsilon,j=1,2.\label{6.010}
\end{eqnarray}
Inserting (\ref{6.010}) into (\ref{6.09}), we have that
\begin{eqnarray}
&&\left\|u_{1}-u_{2}\right\|_{{Z}^{-\frac{1}{2}}([0,\tau])}\leq \frac{1}{2}
\left\|u_{1}-u_{2}\right\|_{{Z}^{-\frac{1}{2}}([0,\tau])}\label{6.011}.
\end{eqnarray}
Thus, $u_{1}=u_{2}$ on $Z^{-\frac{1}{2}}([0,\tau])$. If $\tau=T$, then the claim is proved. If $\tau <T,$ we use the proof of extension.
In particular, $u_{1}(x,y,\tau)=u_{2}(x,y,\tau)$.
Repeating the above proof, we derive that
\begin{eqnarray}
&&\left\|u_{1}-u_{2}\right\|_{{Z}^{-\frac{1}{2}}([\tau,2\tau])}\leq \frac{1}{2}
\left\|u_{1}-u_{2}\right\|_{{Z}^{-\frac{1}{2}}([\tau,2\tau])}\label{6.012}.
\end{eqnarray}
Thus, $u_{1}=u_{2}$ on $Z^{-\frac{1}{2}}([\tau,2\tau])$.If $2\tau=T$, then the claim is proved.
If $2\tau <T,$ we repeatedly use the above proof. We at most repeatedly use $[\frac{T}{\tau}]$ times proof.
We derive that
\begin{eqnarray}
&&\left\|u_{1}-u_{2}\right\|_{{Z}^{-\frac{1}{2}}([[\frac{T}{\tau}]\tau,([\frac{T}{\tau}]+1)\tau])}\leq \frac{1}{2}
\left\|u_{1}-u_{2}\right\|_{{Z}^{-\frac{1}{2}}([[\frac{T}{\tau}]\tau,([\frac{T}{\tau}]+1)\tau])}\label{6.013}.
\end{eqnarray}
Thus, $u_{1}=u_{2}$ on $Z^{-\frac{1}{2}}([[\frac{T}{\tau}]\tau,([\frac{T}{\tau}]+1)\tau]).$
Consequently, $u_{1}=u_{2}$ on $Z^{-\frac{1}{2}}([0,([\frac{T}{\tau}]+1)\tau]).$ Here $([\frac{T}{\tau}]+1)\tau]>T.$

We have completed the proof of Theorem 1.4

\bigskip

\leftline{\large \bf Acknowledgments}

This work is supported by the Natural Science Foundation
of China  under grant numbers 11771127, 11401180, 11571118
and 11471330. The first author is also  supported by
the Young Core Faculty Program of Henan province under
the grant number 5201019430009.
The third author is also  supported by the
Fundamental Research Funds for the Central Universities
of China under the grant number 2017ZD094. We thank
Professors Shuanglin Shao and Hsi-Wei Shih for their
valuable suggestions which considerably improve the
original version of this paper.

\end{document}